\documentclass[a4paper,11pt]{article}

\usepackage{amsmath,amssymb,amsthm}
\usepackage{url,ifthen}
\usepackage{enumitem}
\usepackage[curve,knot]{xypic}

\usepackage[a4paper,left=1in,width=458pt,textheight=650pt,top=1.5in]{geometry}

\usepackage{tikz}
\usetikzlibrary{calc}
\usepgflibrary{shapes.geometric}
\usepgflibrary{shapes.misc}
\usetikzlibrary{positioning}
\usetikzlibrary{decorations}

\usepackage{bbm}
\DeclareMathAlphabet{\mathitbf}{OML}{cmm}{b}{it}

\usepackage{hyperref}

\input cyracc.def

\providecommand{\keywords}[1]{\noindent \footnotesize {\textbf{Keywords:} #1 } \normalsize }

\newboolean{Sectext}
\newsavebox{\MSCtextP}
\newsavebox{\MSCtextS}
\newcommand{\MSC}[2]{
\sbox{\MSCtextP}{#1}
\ifthenelse{\equal{#2}{NONE}}{\setboolean{Sectext}{false}}{\setboolean{Sectext}{true}}
\sbox{\MSCtextS}{#2}
}
\newcommand{\printMSC}{
\ifthenelse{\boolean{Sectext}}{\usebox{\MSCtextP} (Primary), \usebox{\MSCtextS}}{\usebox{\MSCtextP}}
}
\newcommand{\printMSCwithtitle}{
\ifthenelse{\boolean{Sectext}}{\noindent \footnotesize \textbf{Mathematics Subject Classification (2000):} \usebox{\MSCtextP} (Primary), \usebox{\MSCtextS} \normalsize}{\noindent \footnotesize \textbf{Mathematics Subject Classification (2000):} \usebox{\MSCtextP} \normalsize}
}

\hypersetup{
pdftitle={Examples of quantum cluster algebras associated to partial flag varieties},
pdfauthor={Jan E. Grabowski},
pdfstartview=FitH
}

\newcommand{\Ad}{\mathrm{Ad}}
\newcommand{\Adj}[2]{{\Ad_{#1}(#2)}}
\newcommand{\brAdj}[2]{{\underline{\;\;\;\;}\mkern-22mu{{\Ad}_{#1}(#2)}}}
\newcommand{\complex}{\ensuremath \mathbb{C}}
\newcommand{\cross}{\times}
\newcommand{\curly}[1]{\ensuremath{\mathcal{#1}}}
\newcommand{\defeq}{\stackrel{\scriptscriptstyle{\mathrm{def}}}{=}}
\newcommand{\dsum}{\ensuremath{ \oplus}}
\newcommand{\dual}[1]{\ensuremath {#1}^{*}}
\newcommand{\inj}{\hookrightarrow}
\newcommand{\iso}{\ensuremath \cong}
\newcommand{\lgen}{\ensuremath \mathopen{<}} 
\newcommand{\Lie}[1]{\ensuremath{\mathfrak{#1}}}
\newcommand{\nat}{\ensuremath \mathbb{N}}
\newcommand{\qbracket}[3]{\ensuremath{ [\, #1 , #2\, ]_{#3}}}
\newcommand{\qDminor}[1]{D_{q}^{#1}}
\newcommand{\qDminorbf}[1]{\mathitbf{D}_{\mathitbf{q}}^{\mathbf{#1}}}
\newcommand{\qea}[1]{\ensuremath{U_{q}(\Lie{#1})}}
\newcommand{\qebplus}[1]{\ensuremath{U_{q}^{{\scriptscriptstyle \geqslant}}(\Lie{#1})}}
\newcommand{\qinvDminor}[1]{D_{q^{-1}}^{#1}}
\newcommand{\qinvDminorbf}[1]{\mathitbf{D}_{\mathitbf{q}^{\mathbf{-1}}}^{\mathbf{#1}}}
\newcommand{\qminor}[1]{\Delta_{q}^{#1}}
\newcommand{\qminorbf}[1]{\mathbf{\Delta}_{\mathitbf{q}}^{\mathbf{#1}}}
\newcommand{\rgen}{\ensuremath \mathclose{>}}
\newcommand{\union}{\mathrel{\cup}}

\theoremstyle{plain}

\theoremstyle{remark}
\newtheorem{example}{Example}

\title{Examples of quantum cluster algebras \\ associated to partial flag varieties}
\author{Jan E. Grabowski\footnotemark[2] 
\\ \small{\textit{Mathematical Institute, University of Oxford,}}
\\ \small{\textit{24-29 St.\ Giles', Oxford, OX1 3LB, United Kingdom}}
}
\date{27th July 2010}

\begin{document}

\maketitle

\renewcommand{\thefootnote}{\fnsymbol{footnote}}
\footnotetext[2]{Email: \url{jan.grabowski@maths.ox.ac.uk}.  Website: \url{http://people.maths.ox.ac.uk/~grabowsk/}}
\renewcommand{\thefootnote}{\arabic{footnote}}
\setcounter{footnote}{0}

\begin{abstract} 
\noindent We give several explicit examples of quantum cluster algebra structures, as introduced by Berenstein and Zelevinsky, on quantized coordinate rings of partial flag varieties and their associated unipotent radicals.  These structures are shown to be quantizations of the cluster algebra structures found on the corresponding classical objects by Gei\ss, Leclerc and Schr\"{o}er, whose work generalizes that of several other authors.  We also exhibit quantum cluster algebra structures on the quantized enveloping algebras of the Lie algebras of the unipotent radicals.
\end{abstract}

\MSC{\footnotesize 20G42}{\footnotesize 16W35, 17B37}
\printMSCwithtitle

\vfill 

\tableofcontents
\vfill

\pagebreak

\section{Introduction}\label{s:intro}

Cluster algebras were introduced in a series of papers by Fomin and Zelevinsky from 2001 onwards (\cite{FZ-CA1}, \cite{FZ-CA2}, \cite{BFZ-CA3}, \cite{FZ-CA4}), providing a framework for combinatorics associated to dual canonical bases of homogeneous coordinate rings, canonical bases of quantum groups and total positivity for semisimple algebraic groups.  One way of thinking of what it means for a (necessarily) commutative algebra to possess a cluster algebra structure is  that it has a particular form of presentation, with many generators (the cluster variables) but relatively simple relations (exchange relations).

Much work has been done in recent years on the theory of cluster algebras, beginning with the fundamental question of classification.  In the original work of Fomin and Zelevinsky, it was shown that cluster algebras have a Cartan--Killing classification (\cite{FZ-CA2}), with associated concepts of rank, finite and infinite types and roots.   A large number of authors have addressed the combinatorial side, studying the relationships with previously-known combinatorial structures such as generalized associahedra (see for example \cite{Fomin-Reading}).

It has also become common to reformulate the definitions of Fomin and Zelevinsky using quivers and this has highlighted relationships with the representation theory of path algebras and related structures.  A good survey on this material is \cite{Keller-CAsurvey}.  

Significant progress in cluster algebra theory has also come from categorification, notably cluster categories of modules whose tilting theory encodes the cluster combinatorics (\cite{Caldero-Chapoton}, \cite{BMRRT} and many subsequent authors; see \cite{Keller-CAsurvey} for a survey), as well as other forms of categorification (\cite{GLS-PFV}, \cite{Hernandez-Leclerc}).

Examples of cluster algebras include polynomial algebras (of rank 0) and coordinate algebras, e.g. $\complex[SL_{2}]$ (of cluster algebra type $A_{1}$) and $\complex[SL_{4}/N]$, $N$ upper unitriangular matrices (of cluster algebra type $A_{3}$).

Two notable subsequent developments have inspired this work.  The first was the demonstration of cluster algebra structures on homogeneous coordinate rings of partial flag varieties and their associated unipotent radicals by Gei\ss, Leclerc and Schr\"{o}er (\cite{GLS-PFV}), giving a general construction that includes earlier examples by other authors, particularly Scott (\cite{Scott-Grassmannians}).  The second was the definition by Berenstein and Zelevinsky of a quantum cluster algebra (\cite{BZ-QCA}) and a construction conjectured to give quantum cluster algebra structures to quantized coordinate rings of double Bruhat cells.  Among these is the case $\complex_{q}[G^{e,w_{0}}]$, from which one would directly obtain a quantum cluster algebra structure on a subalgebra of $\complex_{q}[N]$, $N$ the maximal unipotent subgroup of the algebraic group $G$.

The natural question posed by considering these together is whether or not one could find \emph{quantum} cluster algebra structures on the \emph{quantized} coordinate rings of partial flag varieties and unipotent radicals.  That is, whether or not one could produce a quantization of the work of Gei\ss, Leclerc and Schr\"{o}er.  By means of the explicit examples given here, we demonstrate that this is a plausible endeavour, deserving of further attention.

The focus of this paper is these examples, rather than general constructions and theorems or applications, though we hope that these will follow shortly.  We begin, however, in Section~\ref{s:objects} with brief descriptions of various objects that one may associate to what we term sub-root data.  By a sub-root datum we mean a pair of root data $(\curly{C},\curly{C}')$ such that $\curly{C}$ has a subsystem of type $\curly{C}'$; one can think of choosing a sub-diagram of a Dynkin diagram, making allowances for graph automorphisms.  Let $I$ be the indexing set (for the simple roots, fundamental weights, etc.) associated to $\curly{C}$ and $J$ that for $\curly{C}'$, as a subset of $I$.

On the geometric side, a sub-root datum gives rise to a standard parabolic subgroup $P_{J}$ of an algebraic group $G$ of type $\curly{C}$, a partial flag variety $G/P_{J}$, an opposite unipotent radical $N_{I\setminus J}^{-}$, their multi-homogeneous coordinate rings and the quantum analogues of these.  On the algebraic side, there are the Lie algebras of the group and subgroups, universal enveloping algebras of these and their quantum analogues.

The last of these, quantized enveloping algebras of Lie algebras of unipotent subgroups, is the point of connection between this paper and our previous work in \cite{BraidedEnvAlgs}.  There, we defined these objects as certain subalgebras of the full quantized enveloping algebra and studied their properties.  In particular, we showed that they are braided Hopf algebras of a special kind, namely Nichols algebras.  These algebras will be denoted $\qea{n_{\mathit{I\setminus J}}^{-}}$.  

In Section~\ref{ss:CAs}, we recall the basic data and constructions appearing in the definition of cluster algebras (of geometric type) and then describe the quantum analogue in Section~\ref{ss:QCAs}.  In that section we also give an example of a quantum cluster algebra that is not associated to a sub-root datum, namely the quantized coordinate algebra $\complex_{q}[SL_{2}]$.  The cluster algebra structure on $\complex[SL_{2}]$ is often used as a simple but effective illustration of the definitions; the quantized version does the same for quantum cluster algebras in an equally pleasing way.  Their cluster algebra type is $A_{1}$ in both cases.

The main body of this paper is contained in Section~\ref{s:examples}, where we give our examples.  We begin in Section~\ref{ss:CPn} with the sub-root datum $(A_{n},A_{n-1})$, when the partial flag variety is just projective space and the unipotent subgroup is an affine space.  The associated (quantized) coordinate algebras are (quantum) symmetric algebras and hence are rank 0 (quantum) cluster algebras.

Next, in Section~\ref{ss:Gr25} we examine the Grassmannian $\mathrm{Gr}(2,5)$, coming from the sub-root datum $(A_{4},A_{1}\dsum A_{2})$ (deleting the node numbered $2$).  We give a quantum cluster algebra structure on $\complex_{q}[\mathrm{Gr}(2,5)]$ coming from quantum Pl\"{u}cker coordinates and the quantum Pl\"{u}cker relations, quantizing the well-known cluster algebra structure on $\complex[\mathrm{Gr}(2,5)]$ first described by Fomin and Zelevinsky in their original papers.  We use this to describe a quantum cluster algebra structure on $\complex_{q}[N_{\{2\}}^{-}]$, the quantized coordinate ring of the opposite big cell.  We also give a bijection which transfers this quantum cluster algebra structure over to $\qea{n_{\{\mathrm{2}\}}^{-}}$, which is in fact isomorphic to $\complex_{q}[N_{\{2\}}^{-}]$.  All of these have associated cluster algebra type $A_{2}$.

Thirdly, in Section~\ref{ss:A4A2}, we identify a quantum cluster algebra structure on $\qea{n_{\{\mathrm{1,2}\}}^{-}} \subseteq \qea{sl_{\mathrm{5}}}$ (deleting nodes $1$ and $2$ from the Dynkin diagram of type $A_{4}$).  This is a more complicated corank 2 example, where the corank is the difference between the rank of $\curly{C}$ and $\curly{C}'$ and is one measure of the complexity of $\qea{n_{\mathit{I\setminus J}}^{-}}$.  This example is of quantum cluster algebra type $A_{3}$.

This example relied on our earlier work (\cite{BraidedEnvAlgs}), which produced information on generating sets for these algebras.  We had implemented the construction in the computer program GAP (\cite{GAP4}), making use of the package QuaGroup (\cite{QuaGroup}), and hence we were able to carry out computer-assisted calculations to verify quasi-commuting of elements and the quantum exchange relations.  

As one would expect, the quantized enveloping algebras $\qea{n^{-}_{\mathit{I\setminus J}}}$ have a close relation with the quantized coordinate rings $\complex_{q}[N_{I\setminus J}^{-}]$, namely that they are graded dual to one another.  However, more can be true: in some cases, these algebras are in fact isomorphic.  It is not known precisely which cases this is true for but when it is, it obviously implies that if a quantum cluster algebra structure exists on the coordinate ring side, one exists on the enveloping algebra side.  However, we have also been able to produce a quantum cluster algebra structure on $\qea{n_{\{\mathrm{1,2}\}}^{-}}$, when it is not known if this is isomorphic to $\complex_{q}[N^{-}_{\{1,2\}}]$ and at present we do not have a quantum cluster algebra structure for the latter.  This suggests that it is reasonable to think that quantized enveloping algebras of Lie algebras of unipotent radicals might all have quantum cluster algebra structures, in contrast to the classical case where non-commutative enveloping algebras cannot be cluster algebras.

We remark that we have demonstrated examples of quantum cluster algebras of type $A_{n}$ for $n=0,1,2,3$, albeit in several different realizations.


Recently, in work with St\'{e}phane Launois (\cite{Gr2nSchubertQCA}) we have extended the example of Section~\ref{ss:Gr25} to the quantized coordinate rings $\mathbb{K}_{q}[\mathrm{Gr}(2,n)]$ for $\mathbb{K}$ an arbitrary field and $n\geq 3$, showing that they are quantum cluster algebras of type $A_{n-3}$.  From this, we have also obtained quantum cluster algebra structures on the quantum Schubert cells of these Grassmannians.  The quantum Schubert cell associated to the partition $(t,s)$ (where $t\geq s$, $t+s\leq 2n-2$) is of quantum cluster algebra type $A_{s-1}$, independent of $t$.  We have also obtained quantum cluster algebra structures on $\mathbb{K}_{q}[\mathrm{Gr}(3,k)]$ for $\mathbb{K}$ a field containing a square root of $q$ and $k\in \{6,7,8\}$, yielding quantum cluster algebras of types $D_{4}$, $E_{6}$ and $E_{8}$.  (These cases exhaust the finite-type quantum cluster algebra structures on quantum Grassmannians.)  We view this as a further step towards achieving the goal of quantizing the work of Gei\ss, Leclerc and Schr\"{o}er, but it is clear that some new methods will be required, likely including a suitable categorification.

Our particular motivation was to give some additional explicit examples of quantum cluster algebras, drawing on the work of Gei\ss, Leclerc and Schr\"{o}er and bringing it together with our own previous work.  On a larger scale, the enterprise of finding quantum cluster algebras is partly aimed at trying to reach a better understanding of the canonical basis of Lusztig and Kashiwara.  Some progress has been made in relating the classical cluster algebras to semi-canonical bases so it is hoped that the quantized versions will play a similar role with respect to canonical bases.  This hope is strengthened by recent work of Nakajima (\cite{Nakajima-QuiverVarClusterAlgs}), who has related these topics using graded quiver varieties, and Lampe (\cite{Lampe}).

\subsection*{Acknowledgements}

I would like to thank a number of individuals and institutions for their help with this project.  In particular, I am grateful to Jeanne Scott for her explanation of her earlier work on cluster algebras, Grassmannians and quantum minors and to Stefan Kolb for several helpful discussions.  I would also like to thank St\'{e}phane Launois for his hospitality at the University of Kent where some of the material here was refined and subsequently extended.

I am grateful to the organisers of the BLOC (Bristol-Leicester-Oxford Colloquium) for an opportunity to speak at their 44th meeting and to all the participants of that meeting for their very helpful comments.  I would also like to acknowledge the provision of facilities by Keble College and the Mathematical Institute in Oxford and the Isaac Newton Institute in Cambridge where parts of this work were carried out.

I would also like to thank the referee for a number of suggested improvements to this work.

\vfill

\section{Objects associated to sub-root data}\label{s:objects}

\subsection{Sub-root data, partial flag varieties and unipotent radicals}\label{ss:SRD+PFVs}

Throughout, we will work over the field $\complex$ and assume that $q$ is not a root of unity.  (It will be clear that more general considerations are possible but we will restrict ourselves to this situation.)

Let $C$ be an $l\cross l$ generalized Cartan matrix with columns indexed by a set $I$.  Let $(H,\Pi,\Pi^{\vee})$ be a minimal realization of $C$, where $H \iso \complex^{2\left| I \right|-\text{rank}(C)}$, $\Pi=\{ \alpha_{i} \mid i \in I \} \subset \dual{H}$ (the simple roots) and $\Pi^{\vee}=\{ h_{i} \mid i \in I \} \subset H$ (the simple coroots).  Then we say $\curly{C}=(C,I,H,\Pi,\Pi^{\vee})$ is a root datum associated to $C$.  (Lusztig (\cite{LusztigBook}) has a more general definition of a root datum but this one will suffice for our purposes here.)

Given two root data $\curly{C}$ and $\curly{C}'=(C',J,H',\Pi',(\Pi')^{\vee})$ and  an injective map $\iota:J \to I$, we will say $\curly{C}'$ is a sub-root datum of $\curly{C}$ if $C'_{ij}=C_{\iota(i)\iota(j)}$ and there exists an injective linear map $s:H' \to H$ such that $s(h_{i}')=h_{\iota(i)}$.  We will generally suppress $\iota$ in our notation and identify $J$ with $\iota(J)\subseteq I$.  If $\curly{C}'$ is a sub-root datum of $\curly{C}$, we will write $\curly{C}'\subseteq \curly{C}$ and set $D=I\setminus J$.  Pictorially, we are choosing a Dynkin diagram and a sub-diagram; the map $\iota$ removes any ambiguities due to graph automorphisms.    

Associated to $\curly{C}' \subseteq \curly{C}$ are several algebraic and geometric objects:

\begin{itemize}
\item If $G=G(\curly{C})$ is a connected semisimple complex algebraic group associated to $\curly{C}$, $G$ has a (standard) parabolic subgroup $P_{J}$ associated to $J\subseteq I$.

\hspace*{1.5em} e.g. $G=G(A_{4})=SL_{5}(\complex)$, $J=\{1,3,4\}\subseteq I=\{1,2,3,4\}$, $D=\{ 2 \}$,
\[ P_{J}=\left\{ \begin{pmatrix} \ast & \ast & \ast & \ast & \ast \\ \ast & \ast & \ast & \ast & \ast \\ 0 & 0 & \ast & \ast & \ast \\ 0 & 0 & \ast & \ast & \ast \\ 0 & 0 & \ast & \ast & \ast \end{pmatrix} \right\}. \]
\item We can form $G/P_{J}$, a partial flag variety.

\hspace*{1.5em} e.g. $G/P_{\emptyset}=G/B$, the full flag variety.  (Here $\curly{C}'=A_{0}$.)

\hspace*{1.5em} e.g. $G(A_{n})/P_{\{1,\ldots ,n-1\}}=\complex \mathbb{P}^{n}$, complex projective space.

\hspace*{1.5em} e.g. $G(A_{4})/P_{\{1,3,4\}}=\mathrm{Gr}(2,5)$, the Grassmannian of planes in $\complex^5$.

The partial flag variety $G/P_{J}$ is a projective variety, via the well-known Pl\"{u}cker embedding $G/P_{J} \inj \prod_{d\in D} \mathbb{P}(L(\omega_{d}))$.  Here, $L(\lambda)$ is the irreducible $G$-module corresponding to a dominant integral weight $\lambda$ and $\{\omega_{i}\}_{i\in I}$ are the fundamental weights.

\item Each standard parabolic subgroup $P_{J}$ has an associated opposite unipotent radical, $N_{D}^{-}$ (see for example \cite{Borel} for details). \label{Nsub2}

\hspace*{1.5em} e.g. $G=G(A_{4})$, $D=\{2 \}$,
\[ N^{-}_{D}=\left\{ \begin{pmatrix} 1 & 0 & 0 & 0 & 0 \\ 0 & 1 & 0 & 0 & 0 \\ \ast & \ast & 1 & 0 & 0 \\ \ast & \ast & 0 & 1 & 0 \\ \ast & \ast & 0 & 0 & 1 \end{pmatrix} \right\}. \]

Note that $N_{D}^{-}$ may be embedded in $G/P_{J}$ as a dense, open subset and is an affine space of the same dimension as $G/P_{J}$.  Indeed, $N_{D}^{-}$ may be identified with the opposite of the big cell of $G/P_{J}$ for its Bruhat decomposition. 

\item Via the Pl\"{u}cker embedding, we may form the corresponding $\nat^{D}$-graded multi-homogeneous coordinate algebras $\complex[N_{D}^{-}]$ and $\complex[G/P_{J}]$, the latter having the decomposition $\bigoplus_{\lambda \in \nat^{D}} \dual{L(\lambda)}$.

\item The coordinate ring $\complex[G]$ has a quantum analogue, $\complex_{q}[G]$ (see for example \cite{Brown-Goodearl}, where this algebra is denoted $\curly{O}_{q}(G)$).  Via this quantized coordinate ring, we can define a quantization of $\complex[G/P_{J}]$, $\complex_{q}[G/P_{J}]$ and also a certain localisation $\complex_{q}[G/P_{J}](f_{e}^{-1})$ whose degree 0 part is the quantized coordinate ring of the corresponding big cell, i.e. the opposite of $N_{D}^{-}$.  We refer the reader to \cite[Section~3]{Kolb} for a survey containing detailed definitions, which will not be required here.  (There, $\complex_{q}[G/P_{J}]$ is denoted $S_{q}[G/P_{S}]$).

\item The algebraic group $G=G(\curly{C})$ has an associated complex Lie algebra, $\Lie{g}=\Lie{g}(\curly{C})$, which has a quantized enveloping algebra $\qea{g}$ defined over $\complex$, with $q\in \dual{\complex}$.   This is generated by elements $F_{i}$, $K_{i}^{\pm 1}$ and $E_{i}$ for $i\in I$, with well-known relations (see for example \cite{Jantzen}).

\item Analogous to the parabolic subgroup $P_{J}$ is the subalgebra $\qea{p_{\mathit{J}}}\defeq \qea{g'}\qebplus{g}\subseteq \qea{g}$, where $\Lie{g}'=\Lie{g}'(\curly{C'})$ and $\qebplus{g}=\lgen K_{i}^{\pm 1}, E_{i} \mid i\in I \rgen$.  In \cite{BraidedEnvAlgs}, we showed that one can associate to this subalgebra a braided Hopf algebra $\qea{n_{\mathit{D}}^{-}}$, dual to $\complex_{q}[N_{D}^{-}]$.  
\end{itemize}

\subsection{Cluster algebras}\label{ss:CAs}

We will recall briefly the definition of a cluster algebra of geometric type with coefficients (\cite{FZ-CA1}).  We start with an \emph{initial seed} $(\underline{y},B)$, consisting of a tuple of generators (called a \emph{cluster}) for the cluster algebra and an \emph{exchange matrix} $B=(b_{ij})$.  (A cluster is not a complete set of generators, but a subset of such a set.)  More seeds are obtained via \emph{mutation} of the initial seed and the elements of the initial cluster are of two types, \emph{mutable cluster variables} and \emph{coefficients}.  The coefficients are not mutated and so are present in every cluster obtained from the initial one by mutation.  An abstract cluster algebra is typically considered to be an algebra over the field of rational functions in the coefficients, taken over a base field $\mathbb{K}$.  We will indicate the mutable variables in a cluster by boldface type.

Mutation has two aspects, matrix mutation and variable mutation.  Matrix mutation $\mu_{k}$ is involutive and given by the rule
\[ (\mu_{k}(B))_{ij} = \begin{cases} -b_{ij} & \text{if}\ i=k\ \text{or}\ j=k \\ b_{ij}+\frac{|b_{ik}|b_{kj}+b_{ik}|b_{kj}|}{2} & \text{otherwise} \end{cases} \]
(Here $k$ is not permitted to be an index of a coefficient variable.)  For example,
\begin{align*}
\begin{pmatrix} 0 & 1 & 0 \\ -1 & 0 & 1 \\ 0 & -1 & 0 \end{pmatrix} & \stackrel{\mu_{1}}{\longrightarrow} \begin{pmatrix} 0 & -1 & 0 \\ 1 & 0 & 1 \\ 0 & -1 & 0 \end{pmatrix} & \begin{pmatrix} 0 & 1 & 0 \\ -1 & 0 & 1 \\ 0 & -1 & 0 \end{pmatrix} & \stackrel{\mu_{2}}{\longrightarrow} \begin{pmatrix} 0 & -1 & 1 \\ 1 & 0 & -1 \\ -1 & 1 & 0 \end{pmatrix}
\end{align*}

\noindent If $(\underline{y}=(y_{1},\ldots ,y_{d}),B)$ is the initial seed then the mutated seed in direction $k$ is given by $(\mu_{k}(\underline{y})=(y_{1},\ldots,y_{k}^{\ast},\ldots,y_{d}),\mu_{k}(B))$, where the new generator $y_{k}^{\ast}$ is determined by the \emph{exchange~relation} 
\[ y_{k}y_{k}^{\ast}=\prod_{b_{ik}>0} y_{i}^{b_{ik}} + \prod_{b_{ik}<0} y_{i}^{-b_{ik}} \]

The alternative quiver description converts $B$ to a quiver by the rule that a strictly positive entry $b_{ij}$ determines a weighted arrow $i\stackrel{b_{ij}}{\to} j$ and a strictly negative one a weighted arrow in the opposite direction.  (Thus $B$ is what is termed a signed adjacency matrix.)  In this approach, coefficients correspond to ``frozen'' vertices, indicated by drawing a box around the vertex.  Then matrix mutation defines the operation of quiver mutation, for example
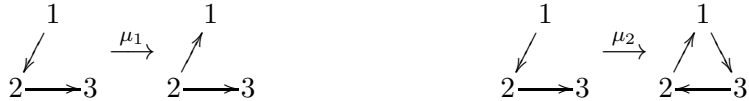
\begin{align*}
\raisebox{1.5em}[3em][2em]{\xymatrix@1@=1pt@!{ {} & {1} \ar[ddl] & {} \\ {} & {} & {} \\ {2} \ar[rr] & {} & {3} }} & \stackrel{\mu_{1}}{\longrightarrow} \raisebox{1.5em}[0em][2em]{\xymatrix@1@=1pt@!{ {} & {1} & {} \\ {} & {} & {} \\ {2} \ar[rr] \ar[uur] & {} & {3} }} & \raisebox{1.5em}[0em][2em]{\xymatrix@1@=1pt@!{ {} & {1} \ar[ddl] & {} \\ {} & {} & {} \\ {2} \ar[rr] & {} & {3} }} & \stackrel{\mu_{2}}{\longrightarrow} \raisebox{1.5em}[0em][2em]{\xymatrix@1@=1pt@!{ {} & {1} \ar[ddr] & {} \\ {} & {} & {} \\ {2} \ar[uur] & {} & {3} \ar[ll] }} 
\end{align*}

We say a $\mathbb{K}$-algebra $\curly{A}$ is a cluster algebra or admits a cluster algebra structure if there is an initial seed such that the set of all cluster variables (i.e. the union of all the clusters) obtained under iterated mutation is a generating set for $\curly{A}$.  In many examples, it is observed that it is not necessary to be able to invert the coefficients so that one can consider the cluster algebra as being defined over the base field $\mathbb{K}$, rather than having to pass to a localisation of $\curly{A}$.

A cluster algebra is of finite type (as all our examples will be) if the quiver of $B$ lies in the same mutation equivalence class as an orientation of a finite-type Dynkin diagram and the type of the cluster algebra is the type of this diagram.  If the cluster algebra under consideration is of finite type $X_{n}$, there is a bijection between the set of all mutable cluster variables (from all clusters) and the almost positive roots of the root system of type $X_{n}$.  (The almost positive roots are the positive roots together with the negative simple roots.)

Work of Gei\ss, Leclerc and Schr\"{o}er (\cite{GLS-PFV}) has identified cluster algebra structures on $\complex[G/P_{J}]$ and $\complex[N_{D}^{-}]$ associated to $\curly{C}'\subseteq \curly{C}$ as above.  Their approach produces the cluster algebra structure from a categorification, specifically from a subcategory of $\Lambda(\curly{C})$-mod, the category of modules of the preprojective algebra of type $\curly{C}$.  The complete rigid modules in this subcategory correspond to the clusters of $\complex[N_{D}^{-}]$ and mutation arises from certain short exact sequences.

Their categorification is used to prove the properties required of a cluster algebra but there is a purely combinatorial recipe to construct the initial seed and exchange matrix (\cite[Section~9.3]{GLS-PFV}).  We will not give this explicitly here but will give the cluster algebra structures it produces alongside our quantized examples, for comparison.

Significantly, Gei\ss, Leclerc and Schr\"{o}er have also shown that monomials in the variables appearing in a single cluster are elements of the dual semicanonical basis of $\complex[N_{I}^{-}]$.  It is conjectured that these monomials lie in the dual \emph{canonical} basis.

\subsection{Quantum cluster algebras}\label{ss:QCAs}

Berenstein and Zelevinsky (\cite{BZ-QCA}) have given a definition of a quantum cluster algebra.  These algebras are now non-commutative but not so far from being commutative.  Each quantum seed includes an additional piece of data, a skew-symmetric matrix $L=(l_{ij})$ describing \emph{quasi-commutation} relations between the variables in the cluster.  Quasi-commuting means $ab=q^{l_{ab}}ba$, also written in $q$-commutator notation as $[a,b]_{q^{l_{ab}}}=0$.  

There is also a mutation rule for these quasi-commutation matrices and a modified exchange relation that involves further coefficients that are powers of $q$ derived from $B$ and $L$, which we describe now.  Given a quantum cluster $\underline{y}=(X_{1},\ldots,X_{r})$, exchange matrix $B$ and quasi-commutation matrix $L$, the exchange relation for mutation in the direction $k$ is given by
\[ X_{k}'=M(-\underline{\boldsymbol{e}}_{k}+\sum_{b_{ik}>0}b_{ik}\underline{\boldsymbol{e}}_{i})+M(-\underline{\boldsymbol{e}}_{k}-\sum_{b_{ik}<0}b_{ik}\underline{\boldsymbol{e}}_{i}) \]
where the vector $\underline{\boldsymbol{e}}_{i}\in \complex^{r}$ ($r$ being the number of rows of $B$) is the $i$th standard basis vector and 
\[ M(a_{1},\dotsc ,a_{r})=q^{\frac{1}{2}\sum_{i<j} a_{i}a_{j}l_{ji}}X_{1}^{a_{1}}\dotsm X_{r}^{a_{r}}. \]
By construction, the integers $a_{i}$ are all non-negative except for $a_{k}=-1$.  The monomial $M$ (as we have defined it here) is related to the concept of a toric frame, also introduced in \cite{BZ-QCA}.  The latter is a technical device used to make the general definition of a quantum cluster algebra.  For our examples, where we start with a known algebra and want to exhibit a quantum cluster algebra structure on this, it will suffice to think of $M$ simply as a rule determining the exchange monomials.  

We note that the presence of the ``$\frac{1}{2}$'' factor in the definition of $M$ suggests that if we wanted to work over fields other than $\complex$, we may need to extend scalars by introducing a square root of $q$.  In fact this will not be necessary in all examples but it would be required in some. 

The natural requirement that all mutated clusters also quasi-commute leads to a compatibility condition between $B$ and $L$, namely that $B^{T}L$ consists of two blocks, one diagonal with positive integer diagonal entries and one zero.  (However, these blocks need not be contiguous, depending on the ordering of the row and column labels.)  We will denote by $0_{m,n}$ the $m\cross n$ zero matrix and by $I_{m}$ the $m\cross m$ identity matrix.  Block matrices will be written in the usual way, e.g. $(A\ B\ C)$.

These definitions are all demonstrated in the following simple yet instructive example.
\begin{example}\label{CqSL2}
Let $\curly{A}=\complex[b,c]$ and define two algebras
\[ \curly{A}_{q}(a,d)=\curly{A}\lgen a,d \rgen / \lgen \qbracket{a}{b}{q},\qbracket{a}{c}{q},\qbracket{d}{b}{q^{-1}},\qbracket{d}{c}{q^{-1}}\rgen\] and
\[ \complex_{q}[SL_{2}]=\curly{A}_{q}(a,d)/\lgen ad=1+qbc,\ da=1+q^{-1}bc \rgen. \]
$\complex_{q}[SL_{2}]$ is a quantum cluster algebra of type $A_{1}$ with initial seed $((\mathitbf{a},b,c),B,L)$,

\begin{align*}
B & = \begin{pmatrix} 0 \\ -1 \\ -1 \end{pmatrix} & & \Gamma(B)=\ \raisebox{2em}[0em][2.5em]{\xymatrix{ {\mathitbf{a}} \ar[r] \ar[d] &*+[F]{b} \\*+[F]{c} & }} & L & = \begin{pmatrix} 0 & 1 & 1 \\ -1 & 0 & 0 \\ -1 & 0 & 0 \end{pmatrix} & \Gamma(L)=\ \raisebox{2em}[0em][0pt]{\xymatrix{ {\mathitbf{a}} \ar[r] \ar[d] &*+[F]{b} \\*+[F]{c} & }}
\end{align*}
Then $B^{T}L=( 2I_{1}\ 0_{1,2})=(2\ 0\ 0)$. The clusters are $(\mathitbf{a},b,c)$ and $(\mathitbf{d},b,c)$.  Setting $X_{1}=\mathitbf{a}$, $X_{2}=b$ and $X_{3}=c$, the first quantum exchange relation is
\begin{eqnarray*}
X_{1}' & = & M(-1,0,0)+M(-1,1,1) \\
 & = & q^{0}X_{1}^{-1}+q^{\frac{1}{2}(-l_{21}-l_{31}+l_{32})}X_{1}^{-1}X_{2}X_{3} \\
 & = & \mathitbf{a}^{-1}+q\mathitbf{a}^{-1}bc,
\end{eqnarray*}
that is, $X_{1}'=\mathitbf{d}$ and $\mathitbf{ad}=1+qbc$.  Note that $\complex[SL_{2}]=\curly{A}_{1}(a,d)/\lgen ad=1+bc \rgen$ has the natural ``classical limit'' cluster algebra structure (forgetting $L$).  

We remark that we have been unable to identify a quantum cluster algebra structure on $\qea{sl_{\mathrm{2}}}$.  The dualisation of the structure we presented here would suggest the mutable cluster variables should be the usual generators $E$ and $F$ but the commutation relation between these is not of a form which fits into the exchange relation type.  However, it is possible that a different presentation might allow $\qea{sl_{\mathrm{2}}}$ to be seen as a quantum cluster algebra.
\end{example}

Importantly, Berenstein and Zelevinsky show that the exchange graph (the graph whose vertices are the clusters and edges are mutations) remains unchanged in the quantum setting.  That is, the matrix $L$ does not influence the exchange graph.  It follows that quantum cluster algebras are classified by Dynkin types in exactly the same way as the classical cluster algebras.

Examples of quantum cluster algebras include quantum symmetric algebras (of rank 0) and it is conjectured that quantum double Bruhat cells are also examples (\cite{BZ-QCA}).  The major aim of this work was to add to this list by exhibiting more examples, which we explain next.

\section{Examples}\label{s:examples}

\subsection{Example: complex projective space}\label{ss:CPn}

The partial flag variety obtained from $G=G(A_{n})=SL_{n+1}(\complex)$, $J=I \setminus \{n\}$ is $G/P_{J}=\complex \mathbb{P}^{n}$, complex projective space.  The corresponding quantized coordinate ring $\complex_{q}[\complex \mathbb{P}^{n}]$ is $S_{q}(\complex^{n+1})$, a quantum symmetric algebra, thus of rank 0 as a quantum cluster algebra.  

The unipotent radical $N_{\{n\}}^{-}$ is $\complex^{n}$, i.e. affine space of dimension $n$, and its quantized coordinate ring $\complex_{q}[\complex^{n}]$ is $S_{q}(\complex^{n})$, so is also a rank 0 quantum cluster algebra.  

The dual to this, $\qea{n_{\{\mathit{n}\}}^{-}}$, is again a quantum symmetric algebra on $n$ variables: the Lie algebra $\Lie{n}^{-}_{\{n\}}$ is the $n$-dimensional natural $\Lie{sl}_{n}(\complex)$-module $V$ with the zero Lie bracket, having universal enveloping algebra $U(\Lie{n}^{-}_{\{n\}})\iso S(V)$.

\subsection{Example: the Grassmannian \texorpdfstring{$\mathrm{Gr}(2,5)$}{Gr(2,5)}}\label{ss:Gr25}

\subsubsection{The partial flag variety}

The partial flag variety obtained from $G=G(A_{4})=SL_{5}(\complex)$, $J=I \setminus \{2\}$ is $G/P_{J}=\mathrm{Gr}(2,5)$, the Grassmannian of 2-dimensional subspaces in $\complex^{5}$.  The quantized coordinate ring $\complex_{q}[\mathrm{Gr}(2,5)]$ is the subalgebra of the quantum matrix algebra $\complex_{q}[M(2,5)]$ generated by the quantum Pl\"{u}cker coordinates.

The quantum matrix algebra is generated by $\{ x_{ij} \mid 1\leq i\leq 2,\ 1\leq j \leq 5 \}$ subject to the quantum $2\cross 2$ matrix relations on each $2\cross 2$ submatrix of \[ \begin{pmatrix} x_{11} & x_{12} & x_{13} & x_{14} & x_{15} \\ x_{21} & x_{22} & x_{23} & x_{24} & x_{25} \end{pmatrix}, \] where the quantum $2\cross 2$ matrix relations on $\left( \begin{smallmatrix} a & b \\ c & d \end{smallmatrix} \right)$ are
\begin{align*} ab & = qba & ac & = qca & bc & = cb \\ bd & = qdb & cd & = qdc & ad-da & = (q-q^{-1})bc. \end{align*}
Hence a presentation for $\complex_{q}[M(2,5)]$ is
\begin{alignat*}{2} \complex_{q}[M(2,5)] = \complex & \lgen x_{ij} & & \mid 1\leq i \leq 2,\ 1\leq j \leq 5 \rgen\  / \\
 & \lgen & & \qbracket{x_{1i}}{x_{1j}}{q}=0,\ \qbracket{x_{2i}}{x_{2j}}{q}=0,\\ & & & \qbracket{x_{1i}}{x_{2i}}{q}=0,\ \qbracket{x_{2i}}{x_{1j}}{\ }=0, \\
 & & & \qbracket{x_{1i}}{x_{2j}}{\ }=(q-q^{-1})x_{1j}x_{2i} \quad \forall\ 1\leq i<j\leq 5 \quad \rgen.
\end{alignat*}

The set of quantum Pl\"{u}cker coordinates $\curly{P}_{q}$ that generate $\complex_{q}[\mathrm{Gr}(2,5)]$ are the $2\cross 2$ quantum minors $\curly{P}_{q} = \{ \Delta_{q}^{ij} \defeq x_{1i}x_{2j}-qx_{1j}x_{2i} \mid 1\leq i < j \leq 5 \}$.  As described in \cite{Leclerc-Zelevinsky} and \cite{Scott-QMinors}, two quantum Pl\"{u}cker coordinates $\Delta_{q}^{ij}$ and $\Delta_{q}^{kl}$ quasi-commute if and only if $\{i,j\}$ and $\{k,l\}$ are weakly separated, meaning---in this particular case---that the corresponding diagonals of a pentagon do not intersect except possibly at vertices.  The power of $q$ appearing in the corresponding quasi-commutation relation is also combinatorially determined.  

We now give an initial quantum seed for a quantum cluster algebra structure on $\complex_{q}[\mathrm{Gr}(2,5)]$.  For the initial quantum cluster we choose
\[ \underline{\tilde{y}} = ( \qminor{15},\qminorbf{14}, \qminorbf{13}, \qminor{12},\qminor{23},\qminor{34},\qminor{45}). \]
This is a set of quasi-commuting variables by the above criterion: the corresponding diagonals of the pentagon are seen to be the five edges (in bijection with the coefficients) and two non-crossing diagonals, $(1,3)$ and $(1,4)$.  That is, this cluster corresponds to a triangulation of the pentagon, as in the classical case (see e.g. \cite{FZ-CA2}).  The elements of this cluster are certainly linearly independent: indeed, the set $\curly{P}_{q}$ of quantum Pl\"{u}cker coordinates is linearly independent (see for example \cite{KLR}).

The corresponding quantum exchange matrix $\tilde{B}$ is equal to that for the well-known cluster algebra structure on $\complex[\mathrm{Gr}(2,5)]$ (\cite{Scott-Grassmannians}) and, along with its quiver $\Gamma(\tilde{B})$, is

\begin{eqnarray*} \tilde{B} = \begin{pmatrix} -1 & 0 \\ \mathbf{0} & \mathbf{-1} \\ \mathbf{1} & \mathbf{0} \\ 0 & 1 \\ 0 & -1 \\ -1 & 1 \\ 1 & 0 \end{pmatrix} &\qquad\qquad  & \Gamma(\tilde{B}) =\ \raisebox{2em}[0em][0em]{\xymatrix{*+[F]{12} \ar[r] & {\mathbf{13}} \ar[r] \ar[d] & {\mathbf{14}} \ar[r] \ar[d] &*+[F]{15} \\ {} &*+[F]{23} &*+[F]{34} \ar[ul] &*+[F]{45} \ar[ul] }}
\end{eqnarray*}
where the quiver vertex corresponding to $\qminor{ij}$ is labelled by $ij$.  We see that this quantum cluster algebra is of type $A_{2}$, since the subquiver on the vertices $\mathbf{13}$ and $\mathbf{14}$ is an orientation of the Dynkin diagram of this type.

The quasi-commutation matrix $\tilde{L}$ and its quiver $\Gamma(\tilde{L})$ are
\begin{eqnarray*} \tilde{L} & = \begin{pmatrix} 0 & -1 & -1 & -1 & 0 & 0 & 1 \\ 1 & \mathbf{0} & \mathbf{-1} & -1 & 0 & 1 & 1 \\ 1 & \mathbf{1} & \mathbf{0} & -1 & 1 & 1 & 2 \\ 1 & 1 & 1 & 0 & 1 & 2 & 2 \\ 0 & 0 & -1 & -1 & 0 & 1 & 2 \\ 0 & -1 & -1 & -2 & -1 & 0 & 1 \\ -1 & -1 & -2 & -2 & -2 & -1 & 0 \end{pmatrix} \\ \Gamma(\tilde{L}) & =\ \raisebox{2em}[6em][7em]
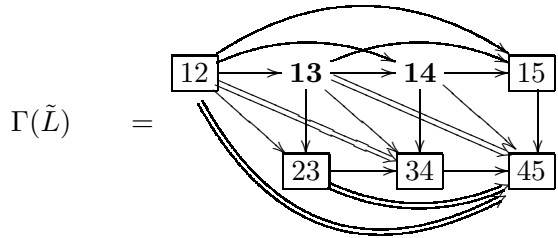
{\xymatrix{*+[F]{12} \ar@<-0.125ex>@/^1pc/[rr] \ar[r] \ar@<0.25ex>@/^2pc/[rrr] \ar@<-1ex>@/_3pc/@{=>}[rrrd] \ar@{=>}[rrd] \ar[rd] & {\mathbf{13}} \ar[r] \ar@<-0.125ex>@/^1pc/[rr] \ar@<0.5ex>@{=>}[rrd] \ar[rd] \ar[d] & {\mathbf{14}} \ar[r] \ar@<0.5ex>[dr] \ar[d] &*+[F]{15} \ar@<0.5ex>[d] \\ {} &*+[F]{23} \ar@<-0.25ex>@/_1pc/@{=>}[rr] \ar[r] &*+[F]{34} \ar[r] &*+[F]{45} }}
\end{eqnarray*}
These are compatible: $\tilde{B}^{T}\tilde{L}=(0_{2,1}\ 2I_{2}\ 0_{2,4})$.

From this data and using the quantum exchange rule, we can write down the exchange relations and identify the remaining cluster variables.  We will give an explicit example of these calculations by considering the mutation $\mu_{2}$, which mutates the minor $\qminorbf{14}$.  (This is $\mu_{2}$, as $\qminorbf{14}$ is the second entry of $\tilde{\underline{y}}$.)  The calculation progresses as in Example~\ref{CqSL2} above:
\begin{align*}
X_{2}' & = M(0,-1,1,0,0,0,1)+M(1,-1,0,0,0,1,0) \\
       & = q^{a}X_{2}^{-1}X_{3}X_{7}+q^{b}X_{1}X_{2}^{-1}X_{6} 
\end{align*}
for some $a,b$.  Now using the relation $X_{1}X_{2}=q^{-1}X_{2}X_{1}$ (which is encoded by $\tilde{L}_{12}$), we may commute the $X_{2}^{-1}$ in the second term to the left and hence multiply through by $X_{2}$ to obtain
\begin{align*}
X_{2}X_{2}' & = q^{a}X_{3}X_{7}+q^{b+1}X_{1}X_{6} \\
           & = q^{-1}X_{3}X_{7}+qX_{1}X_{6}
\end{align*}
since 
\begin{align*}
a & = \tfrac{1}{2}((-1)(1)(\tilde{L}_{32})+(-1)(1)(\tilde{L}_{72})+(1)(1)(\tilde{L}_{73})) = \tfrac{1}{2}(-1+1-2)= -1 \quad\! \text{and} \\
b & = \tfrac{1}{2}((1)(-1)(\tilde{L}_{21})+(1)(1)(\tilde{L}_{61})+(-1)(1)(\tilde{L}_{62})) = \tfrac{1}{2}(-1+0+1)= 0.
\end{align*}
Then, substituting in, we see that $\qminorbf{14}X_{2}'=q^{-1}\qminorbf{13}\qminor{45}+q\qminor{15}\qminor{34}=\qminorbf{14}\qminorbf{35}$, where the latter equality follows from the quantum Pl\"{u}cker relations.  So we deduce that $X_{2}'=\qminorbf{35}$, since $\complex_{q}[\mathrm{Gr}(2,5)]$ is a domain (\cite{KLR}).

We know from the general theory of type $A_{2}$ cluster algebras that only two more cluster variables need to be identified.  These are obtained from the mutations $\mu_{3}$ and $\mu_{2} \circ \mu_{3}$ and the exchange relations determining these are
\begin{align*}
\mu_{3}: & \hspace{0em} & X_{3}X_{3}' & = q^{-1}X_{4}X_{6}+qX_{2}X_{5} & & \implies & X_{3}' & = \qminorbf{24} \\
\mu_{2}\circ \mu_{3}: & \hspace{0em} & X_{2}X_{2}'' & = q^{-1}X_{4}X_{7}+qX_{1}X_{3}' & & \implies & X_{2}'' & = \qminorbf{25}.
\end{align*}
We obtain $X_{3}'=\qminorbf{24}$ and $X_{2}''=\qminorbf{25}$ in the same way as we obtained $X_{2}'$.

In Figure~\ref{Gr25-ex-graph}, we show the exchange graph in this case, with clusters identified with triangulations of the pentagon in the manner described previously and mutations indicated by labelled arrows.  The top vertex corresponds to the initial cluster described here, whose mutable variables are the quantum minors $\qminorbf{14}$ and $\qminorbf{13}$: we number the pentagon's vertices starting with 1 at the top and increasing clockwise.

\begin{figure}
\begin{center} 
\scalebox{0.75}{\begin{tikzpicture}
\node (Pentagon1) at (90:4cm) [regular polygon, regular polygon sides=5, draw,minimum size=2cm] {};
\node (Pentagon2) at (72+90:4cm) [regular polygon, regular polygon sides=5, draw,minimum size=2cm] {};
\node (Pentagon3) at (144+90:4cm) [regular polygon, regular polygon sides=5, draw,minimum size=2cm] {};
\node (Pentagon4) at (216+90:4cm) [regular polygon, regular polygon sides=5, draw,minimum size=2cm] {};
\node (Pentagon5) at (288+90:4cm) [regular polygon, regular polygon sides=5, draw,minimum size=2cm] {};
\draw[->,thick,shorten <=3mm,shorten >=3mm] (Pentagon1.corner 2) -- (Pentagon2.corner 1) node [midway, above] {2};
\draw[->,thick,shorten <=3mm,shorten >=3mm] (Pentagon2.corner 5) -- (Pentagon1.corner 3) node [midway,below] {3};
\draw[->,thick,shorten <=3mm,shorten >=3mm] (Pentagon2.corner 3) -- (Pentagon3.corner 2) node [midway,left] {3};
\draw[->,thick,shorten <=3mm,shorten >=3mm] (Pentagon3.corner 1) -- (Pentagon2.corner 4) node [midway,right] {2};
\draw[->,thick,shorten <=3mm,shorten >=3mm] (Pentagon3.corner 4) -- (Pentagon4.corner 3) node [midway,below] {2};
\draw[<-,thick,shorten <=3mm,shorten >=3mm] (Pentagon3.corner 5) -- (Pentagon4.corner 2) node [midway,above] {3};
\draw[->,thick,shorten <=3mm,shorten >=3mm] (Pentagon4.corner 5) -- (Pentagon5.corner 4) node [midway,right] {3};
\draw[<-,thick,shorten <=3mm,shorten >=3mm] (Pentagon4.corner 1) -- (Pentagon5.corner 3) node [midway,left] {2};
\draw[->,thick,shorten <=3mm,shorten >=3mm] (Pentagon5.corner 1) -- (Pentagon1.corner 5) node [midway,above] {2};
\draw[<-,thick,shorten <=3mm,shorten >=3mm] (Pentagon5.corner 2) -- (Pentagon1.corner 4) node [midway,below] {3};
\draw[] (Pentagon1.corner 1) -- (Pentagon1.corner 4);
\draw[] (Pentagon1.corner 1) -- (Pentagon1.corner 3);
\draw[] (Pentagon2.corner 1) -- (Pentagon2.corner 4);
\draw[] (Pentagon2.corner 2) -- (Pentagon2.corner 4);
\draw[] (Pentagon3.corner 2) -- (Pentagon3.corner 4);
\draw[] (Pentagon3.corner 2) -- (Pentagon3.corner 5);
\draw[] (Pentagon4.corner 2) -- (Pentagon4.corner 5);
\draw[] (Pentagon4.corner 3) -- (Pentagon4.corner 5);
\draw[] (Pentagon5.corner 3) -- (Pentagon5.corner 5);
\draw[] (Pentagon5.corner 1) -- (Pentagon5.corner 3);

\end{tikzpicture}}
\end{center}
\caption{Exchange graph for cluster algebra structure on $\complex[\mathrm{Gr}(2,5)]$ and its quantum analogue.}\label{Gr25-ex-graph}
\end{figure}

In all, ten exchange relations should be considered, corresponding to the ten arrows in Figure~\ref{Gr25-ex-graph}.  It is straightforward (if tedious) to write these down, once one has computed the mutated exchange matrices.  Then one may verify that all the exchange relations are quantum Pl\"{u}cker relations and all cluster variables are quantum minors.  (We originally made use of the computer program Magma (\cite{Magma}) to aid this verification but note that a direct proof, for any $n$, is given in our work with Launois (\cite{Gr2nSchubertQCA}).)

Thus the complete set of cluster variables is
\[ \{ \qminorbf{35}, \qminorbf{25},\qminorbf{24},\qminorbf{14},\qminorbf{13} \} \union \{ \qminor{15},\qminor{12},\qminor{23},\qminor{34},\qminor{45} \} \]
where we have arranged the mutable cluster variables so that they are in bijection with the almost positive roots of the root system of type $A_{2}$ in the order $(\alpha_{1},\alpha_{1}+\alpha_{2},\alpha_{2},-\alpha_{1},-\alpha_{2})$.  We see that this set is equal to $\curly{P}_{q}$, the set of quantum Pl\"{u}cker coordinates which generates $\complex_{q}[\mathrm{Gr}(2,5)]$.  Hence no localisation of the coefficients is needed and $\complex_{q}[\mathrm{Gr}(2,5)]$ is a quantum cluster algebra over $\complex$.

\subsubsection{The unipotent radical}\label{Gr25Unipotent}

The unipotent radical $N_{\{2\}}^{-}$ associated to the above data is an affine space of dimension 6: as a subgroup of $G=SL_{5}(\complex)$ it is that shown on page~\pageref{Nsub2}.

As shown in \cite[Corollaries 1 \& 2]{Heckenberger-Kolb}, the quantized coordinate ring of the unipotent subgroup $N_{\{2\}}^{-}$ is identified with $\complex_{q}[M(2,3)]$.  As is well-known in the classical case and shown in \cite[Proposition~2]{Scott-QMinors} in the quantum case, $\complex_{q}[M(2,3)]$ embeds into $\complex_{q}[\mathrm{Gr}(2,5)]$ (the quantum Stieffel--Pl\"{u}cker correspondence), via an analogue of the map from $M(2,3)$ to $\mathrm{Gr}(2,5)$ that sends a $2\cross 3$ matrix $A$ to the row space of the matrix $\left( \begin{smallmatrix} 0 & 1 \\ -1 & 0 \end{smallmatrix}\right) \dsum A$.

The quantum cluster algebra structure on $\complex_{q}[N^{-}_{\{2\}}]$ may be obtained from that on $\complex_{q}[\mathrm{Gr}(2,5)]$ by noting that as a consequence of this map, the former is generated by the set of quantum minors $\curly{P}_{q}\setminus \{\qminor{12}\}$.  This reflects the construction of $\complex_{q}[N^{-}_{\{2\}}]$ as a quotient of $\complex_{q}[\mathrm{Gr}(2,5)]$.  

Thus an initial quantum seed is given by
\allowdisplaybreaks
\begin{align*} \underline{y} & = ( \qminor{15},\qminorbf{14}, \qminorbf{13}, \qminor{23},\qminor{34},\qminor{45}) & & \\
B & = \begin{pmatrix} -1 & 0 \\ \mathbf{0} & \mathbf{-1} \\ \mathbf{1} & \mathbf{0} \\ 0 & -1 \\ -1 & 1 \\ 1 & 0 \end{pmatrix}  &\quad \Gamma(B) & =\ \raisebox{2em}[4.5em][4em]{\xymatrix{{\mathbf{13}} \ar[r] \ar[d] & {\mathbf{14}} \ar[r] \ar[d] &*+[F]{15} \\*+[F]{23} &*+[F]{34} \ar[ul] &*+[F]{45} \ar[ul] }} \\
L & = \begin{pmatrix} 0 & -1 & -1 & 0 & -1 & 0 \\ 1 & \mathbf{0} & \mathbf{-1} & 0 & 0 & 0 \\ 1 & \mathbf{1} & \mathbf{0} & 1 & 0 & 1 \\ 0 & 0 & -1 & 0 & 0 & 1 \\ 1 & 0 & 0 & 0 & 0 & 1 \\ 0 & 0 & -1 & -1 & -1 & 0 \end{pmatrix} & \Gamma(L) & =\ \raisebox{2em}[0em][0em]
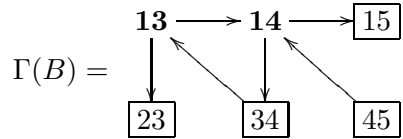
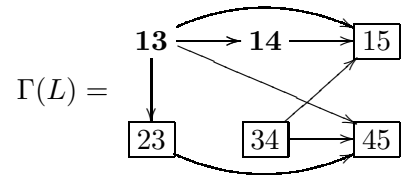
{\xymatrix{{\mathbf{13}} \ar[r] \ar@<0.125ex>@/^1pc/[rr] \ar@<0.5ex>[rrd] \ar[r] \ar[d] & {\mathbf{14}} \ar[r] &*+[F]{15} \\*+[F]{23} \ar@<-0.25ex>@/_1pc/[rr] &*+[F]{34} \ar[r] \ar[ru] &*+[F]{45} }}
\end{align*}
Notice that although $B$ is given simply by deleting the row of $\tilde{B}$ labelled by $12$, the matrix $L$ is not related to $\tilde{L}$ in this way since the quasi-commutation relations are different in the quotient.  The matrices $B$ and $L$ are compatible: $B^{T}L=(0_{2,1}\ 2I_{2}\ 0_{2,3})$.

We can transfer this quantum cluster algebra structure on $\complex_{q}[N_{\{2\}}^{-}]$ to its dual $\qea{n^{-}_{\{\mathrm{2}\}}}$, a subalgebra of $\qea{sl_{\mathrm{5}}}$.  A basis for the degree 1 part of this $\nat$-graded subalgebra, and hence by \cite{BraidedEnvAlgs} a generating set for the whole subalgebra, is as follows, with the action of the generators $E_{i}, F_{i} \in \qea{sl_{\mathrm{2}}\dsum sl_{\mathrm{3}}}$, $i\in \{1,3,4 \}$, indicated:
\[ \xymatrix@1@M=2pt{
{g_{11}\defeq F_{2}K_{2}} \ar@<1ex>[r]^-{F_{1}} \ar@<-1ex>[d]_-{F_{3}} & {g_{21}\defeq \qbracket{F_{1}}{F_{2}}{q}K_{1}K_{2}} \ar@<1ex>[l]^-{E_{1}} \ar@<-1ex>[d]_-{F_{3}} \\
{g_{12}\defeq \qbracket{F_{3}}{F_{2}}{q}K_{2}K_{3}} \ar@<1ex>[r]^-{F_{1}} \ar@<-1ex>[d]_-{F_{4}} \ar@<-1ex>[u]_-{E_{3}} & {g_{22}\defeq \qbracket{F_{3}}{\qbracket{F_{1}}{F_{2}}{q}}{q}K_{1}K_{2}K_{3}} \ar@<1ex>[l]^-{E_{1}} \ar@<-1ex>[d]_-{F_{4}} \ar@<-1ex>[u]_-{E_{3}} \\
{g_{13}\defeq \qbracket{F_{4}}{\qbracket{F_{3}}{F_{2}}{q}}{q}K_{2}K_{3}K_{4}}
 \ar@<1ex>[r]^-{F_{1}} \ar@<-1ex>[u]_-{E_{4}} & {g_{23}\defeq \qbracket{F_{4}}{\qbracket{F_{3}}{\qbracket{F_{1}}{F_{2}}{q}}{q}}{q}K_{1}K_{2}K_{3}K_{4}} \ar@<1ex>[l]^-{E_{1}} \ar@<-1ex>[u]_-{E_{4}}
} \] ($F_{i}$, $K_{i}$ being the usual $\qea{sl_{\mathrm{5}}}$-generators.)  

The relations are
\begin{align*}
\qbracket{g_{1i}}{g_{1j}}{q} = 0 & & & \text{for}\ 1\leq i<j \leq 3 \\
\qbracket{g_{2i}}{g_{2j}}{q} = 0 & & & \text{for}\ 1\leq i<j \leq 3 \\
\qbracket{g_{1i}}{g_{2i}}{q} = 0 & & & \text{for}\ 1\leq i \leq 3 \\
\qbracket{g_{2i}}{g_{1j}}{\ } = 0 & & & \text{for}\ 1\leq i < j \leq 3
\end{align*}
(Note that the last of these is an ordinary commutation relation.)

The quantum cluster algebra structure on $\qea{n^{-}_{\{\mathrm{2}\}}}$ corresponds to that on $\complex_{q}[N^{-}_{\{2\}}]$ via the bijection $\qminor{ij}\mapsto q(\qinvDminor{ij})$, where $\qinvDminor{ij}$ is the $2\cross 2$ $q^{-1}$-minor on columns $i$ and $j$ of the matrix \[ \begin{pmatrix}0 & 1 & g_{11} & g_{12} & g_{13} \\ -1 & 0 & g_{21} & g_{22} & g_{23} \end{pmatrix}. \]  So, the initial quantum cluster is
\[ \underline{y}'=( q\qinvDminor{15},\mathitbf{q}\qinvDminorbf{14}, \mathitbf{q}\qinvDminorbf{13}, q\qinvDminor{23},q\qinvDminor{34},q\qinvDminor{45}) \] and the exchange and quasi-commutation matrices are respectively $B$ and $L$ as above, for the quantum cluster algebra structure on $\complex_{q}[N_{\{2\}}^{-}]$.  In particular, the mutable cluster variables are $\mathitbf{q}\qinvDminorbf{13}=q(q^{-1}g_{11})=F_{2}K_{2}$ and $\mathitbf{q}\qinvDminorbf{14}=g_{12}=\qbracket{F_{3}}{F_{2}}{q}K_{2}K_{3}$.

The remaining\rule{0em}{1em} cluster variables are again obtained from the mutations $\mu_{2}$, $\mu_{3}$ and $\mu_{2} \circ \mu_{3}$ and the three exchange relations determining these are
\begin{align*}
\mu_{2}: & \hspace{0em} & X_{2}X_{2}' & = q^{-1}X_{3}X_{6}+qX_{1}X_{5} & & \implies & X_{2}' & = \mathitbf{q}\qinvDminorbf{35} \\
\mu_{3}: & \hspace{0em} & X_{3}X_{3}' & = X_{5}+qX_{2}X_{4} & & \implies & X_{3}' & = qg_{22}=\mathitbf{q}\qinvDminorbf{24} \\
\mu_{2}\circ \mu_{3}: & \hspace{0em} & X_{2}X_{2}'' & = X_{7}+qX_{1}X_{3}' & & \implies & X_{2}'' & = qg_{23}=\mathitbf{q}\qinvDminorbf{25}.
\end{align*}
where $X_{i}$ is the $i$th element of the initial quantum cluster $\underline{y}'$ above.  Note that the second two of these may be derived from the corresponding relations in $\complex_{q}[\mathrm{Gr}(2,5)]$ by setting $q^{-1}X_{4}=1$ and re-numbering appropriately.  As an example, the first relation, in both the minor notation and in terms of the generators $g_{ij}$, is
\begin{eqnarray*}
\left(\mathitbf{q}\qinvDminorbf{14}\right)\left(\mathitbf{q}\qinvDminorbf{35}\right) & \!\!\! = \!\!\! & q^{-1}\left(\mathitbf{q}\qinvDminorbf{13}\right)\left(q\qinvDminor{45}\right)+q\left(q\qinvDminor{15}\right)\left(q\qinvDminor{34}\right)
 \\ qg_{12}g_{11}g_{23}-g_{12}g_{13}g_{21} &\!\!\! =\!\!\! & g_{11}g_{12}g_{23}-q^{-1}g_{11}g_{13}g_{22}+q^{2}g_{13}g_{11}g_{22}-qg_{13}g_{12}g_{21} \end{eqnarray*}
(an equation of degree 7 in the $F_{i}$).  This and the remaining exchange relations were verified using the package QuaGroup for GAP.

Hence the mutable quantum cluster variables are 
\[ \left\{ \mathitbf{q}\qinvDminorbf{35},\mathitbf{q}\qinvDminorbf{25},\mathitbf{q}\qinvDminorbf{24},\mathitbf{q}\qinvDminorbf{14},\mathitbf{q}\qinvDminorbf{13} \right\} = \left\{ q\qinvDminor{35},qg_{23},qg_{22},g_{12},g_{11} \right\} \]
and the coefficients are
\[ \left\{ q\qinvDminor{15},q\qinvDminor{23},q\qinvDminor{34},q\qinvDminor{45} \right\} = \left\{ g_{13},g_{21},q\qinvDminor{34},q\qinvDminor{45} \right\}
\]
Again the mutable cluster variables are in bijection with the almost-positive roots of $A_{2}$, in the same order as before.  In particular, we see that (scalar multiples of) all generators of $\qea{n^{-}_{\{\mathrm{2}\}}}$ occur in these two lists, so $\qea{n^{-}_{\{\mathrm{2}\}}}$ is a quantum cluster algebra.

\subsection{Example: a corank 2 example in type \texorpdfstring{$A_{4}$}{A4}}\label{ss:A4A2}

The partial flag variety obtained from $G=G(A_{4})=SL_{5}(\complex)$, $J=I \setminus \{1,2\}$ is a projective space of dimension 7.  The unipotent radical $N_{\{1,2\}}^{-}$ associated to the above data is an affine space of dimension 7: as a subgroup of $G=SL_{5}(\complex)$ it is
\[ N^{-}_{\{1,2\}}=\left\{ \begin{pmatrix} 1 & 0 & 0 & 0 & 0 \\ \ast & 1 & 0 & 0 & 0 \\ \ast & \ast & 1 & 0 & 0 \\ \ast & \ast & 0 & 1 & 0 \\ \ast & \ast & 0 & 0 & 1 \end{pmatrix} \right\}. \]

We now describe a quantum cluster algebra structure on $\qea{n^{-}_{\{\mathrm{1,2}\}}}$.  A basis for the degree 1 part of $\qea{n^{-}_{\{\mathrm{1,2}\}}}$ is
\begin{align*}
b_{12} & \defeq F_{1}K_{1} &  b_{24} & \defeq \Adj{F_{3}}{b_{23}}=\qbracket{F_{3}}{F_{2}}{q}K_{2}K_{3} \\ 
b_{23} & \defeq F_{2}K_{2} &  b_{25} & \defeq \Adj{F_{4}}{b_{24}}=\qbracket{F_{4}}{\qbracket{F_{3}}{F_{2}}{q}}{q}K_{2}K_{3}K_{4}
\end{align*}
and these four elements generate $\qea{n^{-}_{\{\mathrm{1,2}\}}}$.  (Their linear independence in $\qea{sl_{\mathrm{5}}}$ follows from the PBW-basis theorem; see for example \cite{Brown-Goodearl}.)

The action of $E_{i},F_{i}\in \qea{sl_{\mathrm{3}}}$ on $\qea{n^{-}_{\{\mathrm{1,2}\}}}_{1}$ may be represented as follows:
\vspace{-0.5em} \[ \xymatrix@1@M=3pt{{b_{12}} & {b_{23}} \ar@<1ex>[r]^{F_{3}} & {b_{24}} \ar@<1ex>[l]^{E_{3}} \ar@<1ex>[r]^{F_{4}} & {b_{25}} \ar@<1ex>[l]^{E_{4}} } \] \vspace{-0.5em}
Note that $\qea{n^{-}_{\{\mathrm{1,2}\}}}_{1}$ is not a simple module.  The relations in $\qea{n^{-}_{\{\mathrm{1,2}\}}}$ are
\begin{align*}
\qbracket{b_{2i}}{b_{2j}}{q} & = 0 & & \text{for}\ 3\leq i<j \leq 5 \\
\brAdj{b_{12}}{\brAdj{b_{12}}{b_{2k}}} & = \qbracket{b_{12}}{\qbracket{b_{12}}{b_{2k}}{q}}{q^{-1}} = 0 & & \text{for}\ 3\leq k \leq 5 \\
\brAdj{b_{2k}}{\brAdj{b_{2k}}{b_{12}}} & = \qbracket{b_{2k}}{\qbracket{b_{2k}}{b_{12}}{q}}{q^{-1}} = 0 & & \text{for}\ 3\leq k \leq 5 
\end{align*}
We note in particular that $\qea{n^{-}_{\{\mathrm{1,2}\}}}$ is not a quadratic algebra, as it has relations in degree 3 for this grading.

Next we give an initial quantum seed $(\underline{y},B,L)$.  For $\underline{y}$ we choose \[ \underline{y} = (\mathitbf{b}_{\mathbf{25}},\mathitbf{b}_{\mathbf{24}},\mathitbf{b}_{\mathbf{23}},\qbracket{b_{12}}{b_{23}}{q},b_{23}b_{12}b_{24}-qb_{24}b_{12}b_{23}, b_{24}b_{12}b_{25}-qb_{25}b_{12}b_{24},\qbracket{b_{25}}{b_{12}}{q}).\]
These elements are linearly independent, as a consequence of the PBW-basis theorem and consideration of their multi-degrees as (homogeneous) polynomials in the $F_{i}$.  In order to interpret these expressions, we set
\begin{align*} b_{13} & = \qbracket{b_{23}}{b_{12}}{q}, & b_{14} & = \qbracket{b_{24}}{b_{12}}{q}, & b_{15} & = \qbracket{b_{25}}{b_{12}}{q}. \end{align*}
Then the elements of $\underline{y}$ can be expressed as either entries or $2\cross 2$ quantum minors of the matrix \[ \begin{pmatrix} 1 & (1-q^{2})b_{12} & b_{13} & b_{14} & b_{15} \\ 0 & 1 & b_{23} & b_{24} & b_{25} \end{pmatrix}.\] Writing $D_{q}^{ij}$ for the $2\cross 2$ quantum minor on columns $i$ and $j$ of the above matrix, we have the following alternative description of $\underline{y}$:
\[ \underline{y}=(\mathitbf{b}_{\mathbf{25}},\mathitbf{b}_{\mathbf{24}},\mathitbf{b}_{\mathbf{23}},\qDminor{23}, \qDminor{34}, \qDminor{45},b_{15}). \]
This quantum cluster is dual to a quantization of a classical cluster for $\complex[N_{\{1,2\}}^{-}]$ constructed by the method of \cite{GLS-PFV}, which is
\[ \underline{z} = (\mathitbf{n}_{\mathbf{25}},\mathitbf{n}_{\mathbf{24}},\mathitbf{n}_{\mathbf{23}},D^{23},D^{34},D^{45},n_{15}) \]
where $D^{ij}$ denotes the (ordinary) $2\cross 2$ minor on columns $i$ and $j$ of the matrix
\[ \begin{pmatrix} 1 & n_{12} & n_{13} & n_{14} & n_{15} \\ 0 & 1 & n_{23} & n_{24} & n_{25} \end{pmatrix}.\]
We could use the minor labels instead, except that $b_{15}$ is not a $2\cross 2$ minor of the above matrix.  It is however a minor of the matrix
\[ \begin{pmatrix} 1 & (1-q^{2})b_{12} & b_{13} & b_{14} & b_{15} & 0 \\ 0 & 1 & b_{23} & b_{24} & b_{25} & 1 \end{pmatrix} \]
and all other relevant minors are unchanged by considering this matrix instead.  Then $b_{15}$ is given by the quantum minor $\qDminor{56}$, leading to our final expression for $\underline{y}$,
\[ \underline{y}=(\qDminorbf{15},\qDminorbf{14},\qDminorbf{13},\qDminor{23},\qDminor{34},\qDminor{45},\qDminor{56}) \]
which gives us the labelling set we use below.

The exchange matrix and quasi-commutation matrix are
\begin{eqnarray*} B = \begin{pmatrix} \mathbf{0} & \mathbf{-1} & \mathbf{0} \\ \mathbf{1} & \mathbf{0} & \mathbf{-1} \\ \mathbf{0} & \mathbf{1} & \mathbf{0} \\ 0 & 0 & -1 \\ 0 & -1 & 1 \\ -1 & 1 & 0 \\ 1 & 0 & 0 \end{pmatrix} &\qquad\qquad  &   L=\begin{pmatrix} \mathbf{0} & \mathbf{-1} & \mathbf{-1} & 0 & -1 & 0 & -1 \\ \mathbf{1} & \mathbf{0} & \mathbf{-1} & 0 & 0 & 0 & 0 \\ \mathbf{1} & \mathbf{1} & \mathbf{0} & 1 & 0 & 1 & 0 \\ 0 & 0 & -1 & 0 & 0 & 1 & 1 \\ 1 & 0 & 0 & 0 & 0 & 1 & 1 \\ 0 & 0 & -1 & -1 & -1 & 0 & 0 \\ 1 & 0 & 0 & -1 & -1 & 0 & 0 \end{pmatrix}
\end{eqnarray*}
We see that this example is of cluster algebra type $A_{3}$; the corresponding quivers are
\vspace{1.5em}
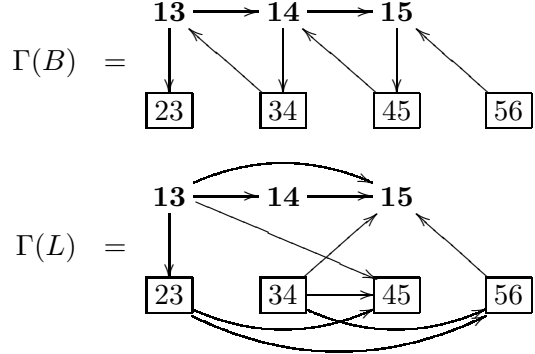
\begin{eqnarray*}
 \Gamma(B) & =\ \raisebox{2em}[0em][3em]{\xymatrix{{\mathbf{13}} \ar[r] \ar[d] & {\mathbf{14}} \ar[r] \ar[d] & {\mathbf{15}} \ar[d]  &{} \\*+[F]{23}  &*+[F]{34} \ar[ul] &*+[F]{45} \ar[ul] &*+[F]{56} \ar[ul] }} \\ \Gamma(L) & =\ \raisebox{2em}[3em][3em]{\xymatrix{{\mathbf{13}} \ar[r] \ar@<0.125ex>@/^1pc/[rr] \ar@<0.5ex>[rrd] \ar[d] & {\mathbf{14}} \ar[r] & {\mathbf{15}} & \\*+[F]{23} \ar@<-0.25ex>@/_1pc/[rr] \ar@<-0.75ex>@/_1.5pc/[rrr] &*+[F]{34} \ar[r] \ar[ru] \ar@<-0.25ex>@/_1pc/[rr] &*+[F]{45} &*+[F]{56} \ar[ul] }}
\end{eqnarray*}
The matrices $B$ and $L$ are compatible: $B^{T}L=(2I_{3}\ 0_{3,4})$.

The remaining cluster variables are obtained from the following set of mutations and corresponding exchange relations.
\begin{align*}
\mu_{1}: & \hspace{0em} & X_{1}X_{1}' & = q^{-1}X_{2}X_{7}+X_{6} & & \implies & X_{1}' & = b_{14}=\qDminorbf{46} \\
\mu_{2}: & \hspace{0em} & X_{2}X_{2}' & = q^{-1}X_{3}X_{6}+qX_{1}X_{5} & & \implies & X_{2}' & = \qDminorbf{35} \\
\mu_{3}: & \hspace{0em} & X_{3}X_{3}' & = X_{5}+qX_{2}X_{4} & & \implies & X_{3}' & = \qDminorbf{24} \\
\mu_{2}\circ \mu_{1}: & \hspace{0em} & X_{2}X_{2}'' & = q^{-1}X_{1}'X_{3}+X_{5} & & \implies & X_{2}'' & = b_{13}=\qDminorbf{36} \\
\mu_{3}\circ \mu_{2}: & \hspace{0em} & X_{3}X_{3}'' & = X_{2}'+qX_{1}X_{4} & & \implies & X_{3}'' & = \qDminorbf{25} \\
\mu_{3}\circ \mu_{2}\circ \mu_{1}: & \hspace{0em} & X_{3}X_{3}''' & = q^{-1/2}X_{2}''+q^{1/2}X_{4} & & \implies & X_{3}''' & = q^{-1/2}(1-q^{2})b_{12} \\ & & & & & & & =\mathitbf{q}^{\mathbf{-1/2}}\qDminorbf{26}
\end{align*}
As usual, $X_{i}$ denotes the $i$th entry of $\underline{y}$.  These exchange relations were again verified using the package QuaGroup for GAP.

We note in particular the appearance of $q^{1/2}$ in the last exchange relation.  This also occurs in Berenstein and Zelevinsky's study of quantum cluster algebra structures for quantum double Bruhat cells (\cite{BZ-QCA}).  They work over $\mathbb{Q}(q)$ and note that the extension of scalars to ensure a square root of $q$ appears to be necessary.  Since we are working over $\complex$, we do not need to make any field extension, of course, but this must be borne in mind if one were to consider other fields.  

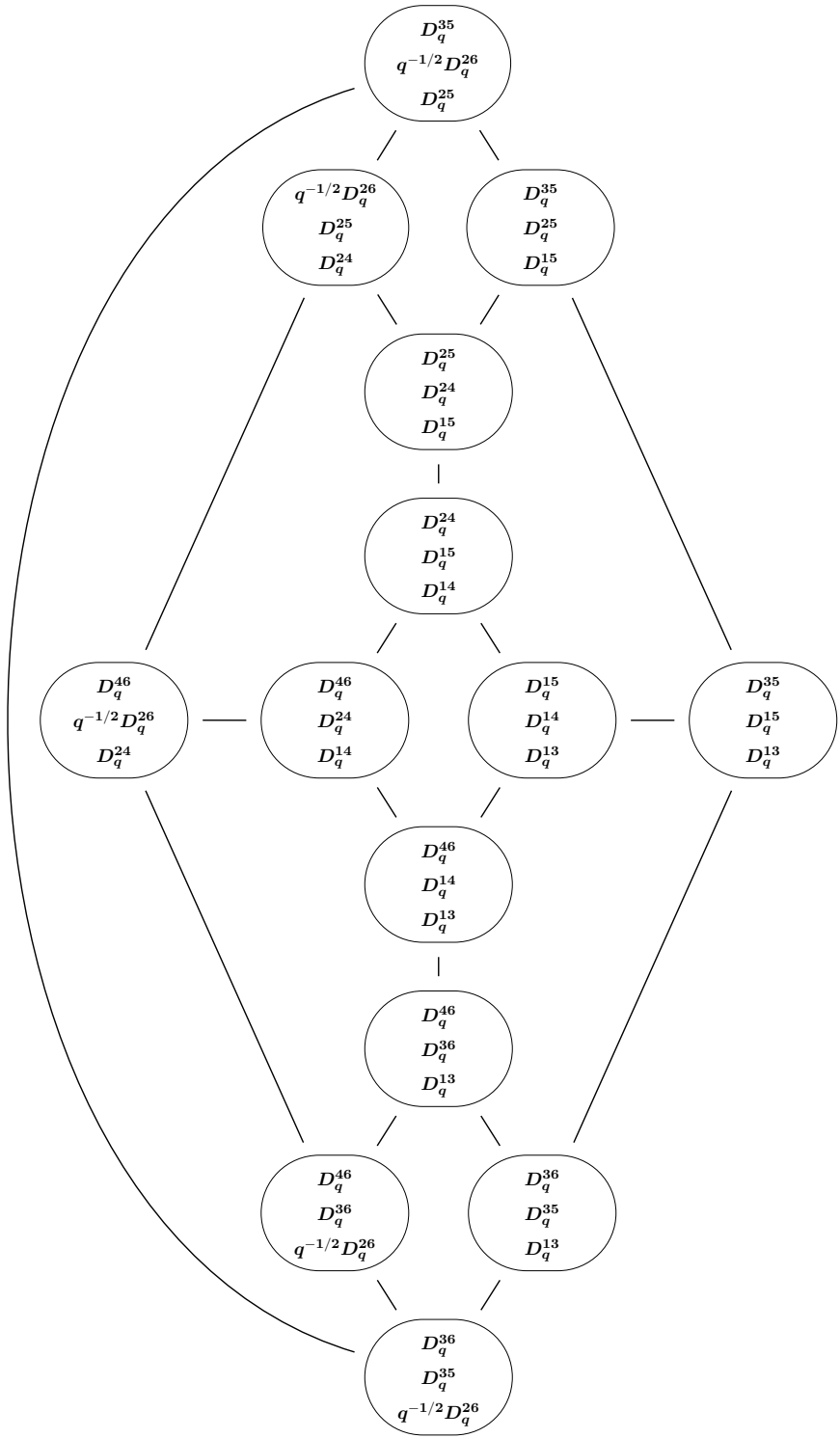
\begin{figure}
\begin{center} 
\renewcommand{\arraystretch}{1.5}
\hspace{-12.5em}\scalebox{0.66}{\begin{tikzpicture}

\begin{scope}[node distance=1cm and 1.5cm]
\node (Hexagon1) [rounded rectangle,draw] at (90:0cm) {$\begin{array}{c} \qDminorbf{35} \\ \mathitbf{q}^{\mathbf{-1/2}}\qDminorbf{26} \\ \qDminorbf{25} \end{array}$};
\node (Hexagon2) [below left=of Hexagon1,rounded rectangle,draw] {$ \begin{array}{c} \mathitbf{q}^{\mathbf{-1/2}}\qDminorbf{26} \\ \qDminorbf{25} \\ \qDminorbf{24}\end{array}$};
\node (Hexagon3) [below right=of Hexagon2,rounded rectangle,draw] {$ \begin{array}{c} \rule{1.3em}{0em} \qDminorbf{25} \rule{1.3em}{0em} \\ \qDminorbf{24}\\ \qDminorbf{15}\end{array}$};
\node (Hexagon4) [below right=of Hexagon1,rounded rectangle,draw] {$ \begin{array}{c} \rule{1.3em}{0em} \qDminorbf{35} \rule{1.3em}{0em} \\ \qDminorbf{25}\\ \qDminorbf{15}\end{array}$};
\node (Hexagon5) [below=of Hexagon3,rounded rectangle,draw] {$ \begin{array}{c} \rule{1.3em}{0em} \qDminorbf{24} \rule{1.3em}{0em} \\ \qDminorbf{15}\\ \qDminorbf{14}\end{array}$};
\node (Hexagon6) [below left=of Hexagon5,rounded rectangle,draw] {$ \begin{array}{c} \rule{1.3em}{0em} \qDminorbf{46} \rule{1.3em}{0em} \\ \qDminorbf{24}\\ \qDminorbf{14}\end{array}$};
\node (Hexagon7) [below right=of Hexagon6,rounded rectangle,draw] {$ \begin{array}{c} \rule{1.3em}{0em} \qDminorbf{46} \rule{1.3em}{0em} \\ \qDminorbf{14}\\ \qDminorbf{13}\end{array}$};
\node (Hexagon8) [below right=of Hexagon5,rounded rectangle,draw] {$ \begin{array}{c} \rule{1.3em}{0em} \qDminorbf{15} \rule{1.3em}{0em} \\ \qDminorbf{14}\\ \qDminorbf{13}\end{array}$};
\node (Hexagon9) [below=of Hexagon7,rounded rectangle,draw] {$ \begin{array}{c} \rule{1.3em}{0em} \qDminorbf{46} \rule{1.3em}{0em} \\ \qDminorbf{36}\\ \qDminorbf{13}\end{array}$};
\node (Hexagon10) [below left=of Hexagon9,rounded rectangle,draw] {$ \begin{array}{c} \rule{1.3em}{0em} \qDminorbf{46} \rule{1.3em}{0em} \\ \qDminorbf{36}\\ \mathitbf{q}^{\mathbf{-1/2}}\qDminorbf{26}\end{array}$};
\node (Hexagon11) [below right=of Hexagon10,rounded rectangle,draw] {$ \begin{array}{c} \rule{1.3em}{0em} \qDminorbf{36} \rule{1.3em}{0em} \\ \qDminorbf{35}\\ \mathitbf{q}^{\mathbf{-1/2}}\qDminorbf{26}\end{array}$};
\node (Hexagon12) [below right=of Hexagon9,rounded rectangle,draw] {$ \begin{array}{c} \rule{1.3em}{0em} \qDminorbf{36} \rule{1.3em}{0em} \\ \qDminorbf{35}\\ \qDminorbf{13}\end{array}$};
\node (Hexagon13) [left=of Hexagon6,rounded rectangle,draw] {$ \begin{array}{c} \rule{1.3em}{0em} \qDminorbf{46} \rule{1.3em}{0em} \\ \mathitbf{q}^{\mathbf{-1/2}}\qDminorbf{26}\\ \qDminorbf{24}\end{array}$};
\node (Hexagon14) [right=of Hexagon8,rounded rectangle,draw] {$ \begin{array}{c} \rule{1.3em}{0em} \qDminorbf{35} \rule{1.3em}{0em} \\ \qDminorbf{15}\\ \qDminorbf{13}\end{array}$};
\end{scope}

\begin{scope}[node distance=7cm and 7.5cm]
\node (ctrl1) [left=of Hexagon2] {};
\node (ctrl2) [left=of Hexagon13] {};
\node (ctrl3) [left=of Hexagon10] {};
\end{scope}

\draw[thick,shorten <=3mm,shorten >=3mm] (Hexagon1) -- (Hexagon2);
\draw[thick,shorten <=3mm,shorten >=3mm] (Hexagon1) -- (Hexagon4);
\draw[thick,shorten <=3mm,shorten >=3mm] (Hexagon2) -- (Hexagon3);
\draw[thick,shorten <=3mm,shorten >=3mm] (Hexagon2) -- (Hexagon13);
\draw[thick,shorten <=3mm,shorten >=3mm] (Hexagon3) -- (Hexagon4);
\draw[thick,shorten <=3mm,shorten >=3mm] (Hexagon3) -- (Hexagon5);
\draw[thick,shorten <=3mm,shorten >=3mm] (Hexagon4) -- (Hexagon14);
\draw[thick,shorten <=3mm,shorten >=3mm] (Hexagon5) -- (Hexagon6);
\draw[thick,shorten <=3mm,shorten >=3mm] (Hexagon5) -- (Hexagon8);
\draw[thick,shorten <=3mm,shorten >=3mm] (Hexagon6) -- (Hexagon7);
\draw[thick,shorten <=3mm,shorten >=3mm] (Hexagon6) -- (Hexagon13);
\draw[thick,shorten <=3mm,shorten >=3mm] (Hexagon7) -- (Hexagon8);
\draw[thick,shorten <=3mm,shorten >=3mm] (Hexagon7) -- (Hexagon9);
\draw[thick,shorten <=3mm,shorten >=3mm] (Hexagon8) -- (Hexagon14);
\draw[thick,shorten <=3mm,shorten >=3mm] (Hexagon9) -- (Hexagon10);
\draw[thick,shorten <=3mm,shorten >=3mm] (Hexagon9) -- (Hexagon12);
\draw[thick,shorten <=3mm,shorten >=3mm] (Hexagon10) -- (Hexagon11);
\draw[thick,shorten <=3mm,shorten >=3mm] (Hexagon10) -- (Hexagon13);
\draw[thick,shorten <=3mm,shorten >=3mm] (Hexagon11) -- (Hexagon12);
\draw[thick,shorten <=3mm,shorten >=3mm] (Hexagon11) .. controls (ctrl3) and (ctrl1) .. (Hexagon1);
\draw[thick,shorten <=3mm,shorten >=3mm] (Hexagon12) -- (Hexagon14);

\end{tikzpicture}}
\renewcommand{\arraystretch}{1}
\end{center}
\caption{Exchange graph for quantum cluster algebra structure on $\qea{n^{-}_{\{{\mathrm{1,2}\}}}}$.}\label{Uqn12-ex-graph}
\end{figure}

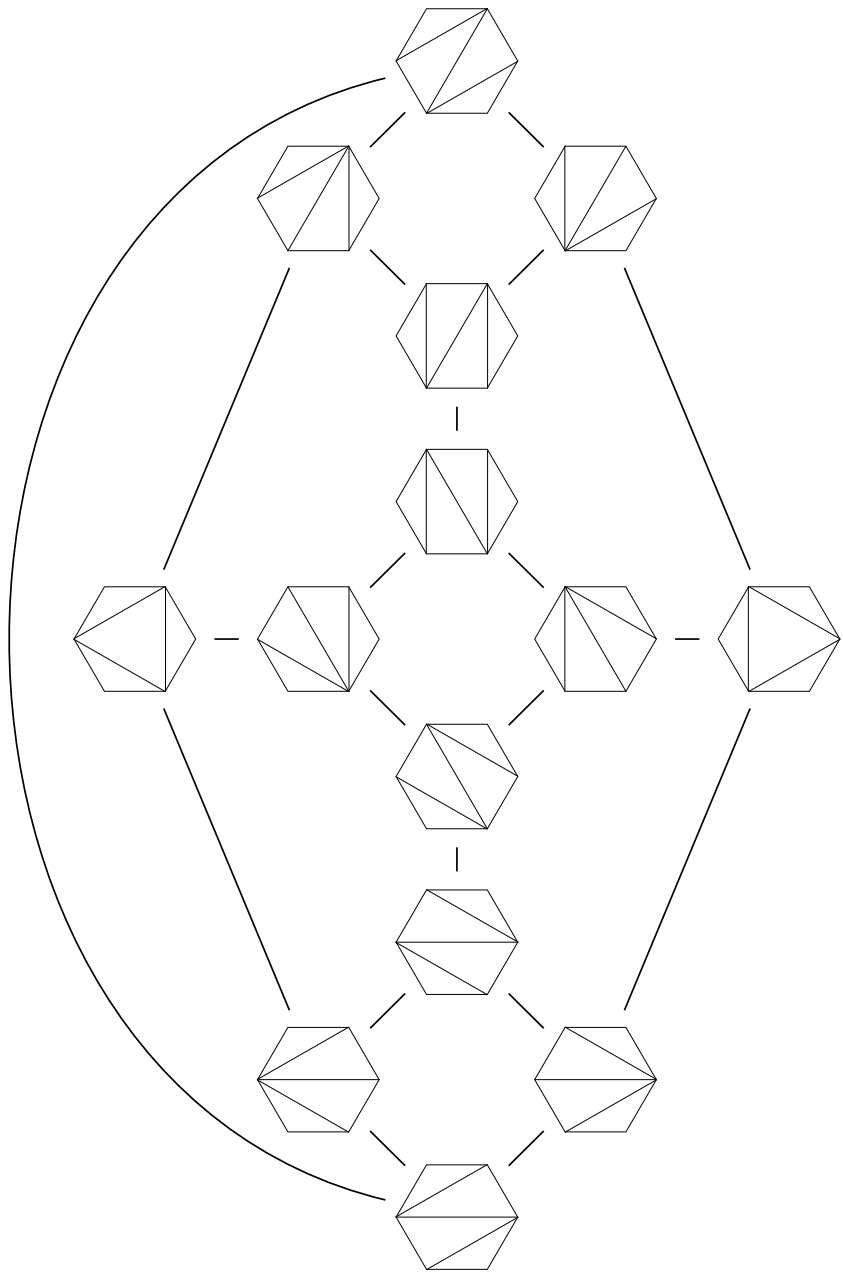
\begin{figure}
\begin{center} 
\hspace{-10.5em}\scalebox{0.8}{\begin{tikzpicture}

\node (Hexagon1) at (90:0cm) [regular polygon, regular polygon sides=6, draw,minimum size=2cm] {};
\node (Hexagon2) [below left=of Hexagon1,regular polygon, regular polygon sides=6, draw,minimum size=2cm] {};
\node (Hexagon3) [below right=of Hexagon2,regular polygon, regular polygon sides=6, draw,minimum size=2cm] {};
\node (Hexagon4) [below right=of Hexagon1,regular polygon, regular polygon sides=6, draw,minimum size=2cm] {};
\node (Hexagon5) [below=of Hexagon3,regular polygon, regular polygon sides=6, draw,minimum size=2cm] {};
\node (Hexagon6) [below left=of Hexagon5,regular polygon, regular polygon sides=6, draw,minimum size=2cm] {};
\node (Hexagon7) [below right=of Hexagon6,regular polygon, regular polygon sides=6, draw,minimum size=2cm] {};
\node (Hexagon8) [below right=of Hexagon5,regular polygon, regular polygon sides=6, draw,minimum size=2cm] {};
\node (Hexagon9) [below=of Hexagon7,regular polygon, regular polygon sides=6, draw,minimum size=2cm] {};
\node (Hexagon10) [below left=of Hexagon9,regular polygon, regular polygon sides=6, draw,minimum size=2cm] {};
\node (Hexagon11) [below right=of Hexagon10,regular polygon, regular polygon sides=6, draw,minimum size=2cm] {};
\node (Hexagon12) [below right=of Hexagon9,regular polygon, regular polygon sides=6, draw,minimum size=2cm] {};
\node (Hexagon13) [left=of Hexagon6,regular polygon, regular polygon sides=6, draw,minimum size=2cm] {};
\node (Hexagon14) [right=of Hexagon8,regular polygon, regular polygon sides=6, draw,minimum size=2cm] {};

\begin{scope}[node distance=6cm and 6cm]
\node (ctrl1) [left=of Hexagon2] {};
\node (ctrl2) [left=of Hexagon13] {};
\node (ctrl3) [left=of Hexagon10] {};
\end{scope}

\draw[thick,shorten <=3mm,shorten >=3mm] (Hexagon1) -- (Hexagon2);
\draw[thick,shorten <=3mm,shorten >=3mm] (Hexagon1) -- (Hexagon4);
\draw[thick,shorten <=3mm,shorten >=3mm] (Hexagon2) -- (Hexagon3);
\draw[thick,shorten <=3mm,shorten >=3mm] (Hexagon2) -- (Hexagon13);
\draw[thick,shorten <=3mm,shorten >=3mm] (Hexagon3) -- (Hexagon4);
\draw[thick,shorten <=3mm,shorten >=3mm] (Hexagon3) -- (Hexagon5);
\draw[thick,shorten <=3mm,shorten >=3mm] (Hexagon4) -- (Hexagon14);
\draw[thick,shorten <=3mm,shorten >=3mm] (Hexagon5) -- (Hexagon6);
\draw[thick,shorten <=3mm,shorten >=3mm] (Hexagon5) -- (Hexagon8);
\draw[thick,shorten <=3mm,shorten >=3mm] (Hexagon6) -- (Hexagon7);
\draw[thick,shorten <=3mm,shorten >=3mm] (Hexagon6) -- (Hexagon13);
\draw[thick,shorten <=3mm,shorten >=3mm] (Hexagon7) -- (Hexagon8);
\draw[thick,shorten <=3mm,shorten >=3mm] (Hexagon7) -- (Hexagon9);
\draw[thick,shorten <=3mm,shorten >=3mm] (Hexagon8) -- (Hexagon14);
\draw[thick,shorten <=3mm,shorten >=3mm] (Hexagon9) -- (Hexagon10);
\draw[thick,shorten <=3mm,shorten >=3mm] (Hexagon9) -- (Hexagon12);
\draw[thick,shorten <=3mm,shorten >=3mm] (Hexagon10) -- (Hexagon11);
\draw[thick,shorten <=3mm,shorten >=3mm] (Hexagon10) -- (Hexagon13);
\draw[thick,shorten <=3mm,shorten >=3mm] (Hexagon11) -- (Hexagon12);
\draw[thick,shorten <=3mm,shorten >=3mm] (Hexagon11) .. controls (ctrl3) and (ctrl1) .. (Hexagon1);
\draw[thick,shorten <=3mm,shorten >=3mm] (Hexagon12) -- (Hexagon14);

\draw[] (Hexagon1.corner 4) -- (Hexagon1.corner 6);
\draw[] (Hexagon1.corner 1) -- (Hexagon1.corner 3);
\draw[] (Hexagon1.corner 1) -- (Hexagon1.corner 4);

\draw[] (Hexagon2.corner 1) -- (Hexagon2.corner 3);
\draw[] (Hexagon2.corner 1) -- (Hexagon2.corner 4);
\draw[] (Hexagon2.corner 1) -- (Hexagon2.corner 5);

\draw[] (Hexagon3.corner 1) -- (Hexagon3.corner 4);
\draw[] (Hexagon3.corner 2) -- (Hexagon3.corner 4);
\draw[] (Hexagon3.corner 1) -- (Hexagon3.corner 5);

\draw[] (Hexagon4.corner 1) -- (Hexagon4.corner 4);
\draw[] (Hexagon4.corner 2) -- (Hexagon4.corner 4);
\draw[] (Hexagon4.corner 4) -- (Hexagon4.corner 6);

\draw[] (Hexagon5.corner 1) -- (Hexagon5.corner 5);
\draw[] (Hexagon5.corner 2) -- (Hexagon5.corner 4);
\draw[] (Hexagon5.corner 2) -- (Hexagon5.corner 5);

\draw[] (Hexagon6.corner 1) -- (Hexagon6.corner 5);
\draw[] (Hexagon6.corner 2) -- (Hexagon6.corner 5);
\draw[] (Hexagon6.corner 3) -- (Hexagon6.corner 5);

\draw[] (Hexagon7.corner 2) -- (Hexagon7.corner 6);
\draw[] (Hexagon7.corner 2) -- (Hexagon7.corner 5);
\draw[] (Hexagon7.corner 3) -- (Hexagon7.corner 5);

\draw[] (Hexagon8.corner 2) -- (Hexagon8.corner 4);
\draw[] (Hexagon8.corner 2) -- (Hexagon8.corner 5);
\draw[] (Hexagon8.corner 2) -- (Hexagon8.corner 6);

\draw[] (Hexagon9.corner 2) -- (Hexagon9.corner 6);
\draw[] (Hexagon9.corner 3) -- (Hexagon9.corner 5);
\draw[] (Hexagon9.corner 3) -- (Hexagon9.corner 6);

\draw[] (Hexagon10.corner 1) -- (Hexagon10.corner 3);
\draw[] (Hexagon10.corner 3) -- (Hexagon10.corner 6);
\draw[] (Hexagon10.corner 3) -- (Hexagon10.corner 5);

\draw[] (Hexagon11.corner 1) -- (Hexagon11.corner 3);
\draw[] (Hexagon11.corner 3) -- (Hexagon11.corner 6);
\draw[] (Hexagon11.corner 4) -- (Hexagon11.corner 6);

\draw[] (Hexagon12.corner 2) -- (Hexagon12.corner 6);
\draw[] (Hexagon12.corner 3) -- (Hexagon12.corner 6);
\draw[] (Hexagon12.corner 4) -- (Hexagon12.corner 6);

\draw[] (Hexagon13.corner 1) -- (Hexagon13.corner 3);
\draw[] (Hexagon13.corner 3) -- (Hexagon13.corner 5);
\draw[] (Hexagon13.corner 5) -- (Hexagon13.corner 1);

\draw[] (Hexagon14.corner 2) -- (Hexagon14.corner 6);
\draw[] (Hexagon14.corner 4) -- (Hexagon14.corner 6);
\draw[] (Hexagon14.corner 2) -- (Hexagon14.corner 4);

\end{tikzpicture}}
\end{center}
\caption{Exchange graph for quantum cluster algebra structure on $\qea{n^{-}_{\{{\mathrm{1,2}\}}}}$, with hexagon triangulations.  (We label the hexagon vertices with 1 at the top left vertex and increasing in a clockwise direction.)}\label{hexagons}
\end{figure}

The complete set of quantum cluster variables is
\[ \{ \qDminorbf{46},\qDminorbf{36},\qDminorbf{35},\mathitbf{q}^{\mathbf{-1/2}}\qDminorbf{26},\qDminorbf{25},\qDminorbf{24},\qDminorbf{15},\qDminorbf{14},\qDminorbf{13} \} \union \{\qDminor{23},\qDminor{34},\qDminor{45},\qDminor{56} \}. \]
The mutable cluster variables are in bijection with the almost-positive roots of the root system of type $A_{3}$, in the following order:
\[ (\alpha_{1},\alpha_{1}+\alpha_{2},\alpha_{2},\alpha_{1}+\alpha_{2}+\alpha_{3},\alpha_{2}+\alpha_{3},\alpha_{3},-\alpha_{1},-\alpha_{2},-\alpha_{3}). \]
In Figure~\ref{Uqn12-ex-graph}, we give the exchange graph for this quantum cluster algebra structure, omitting the coefficients (which of course occur in each cluster).  The exchange graph for a cluster algebra of type $A_{3}$ is isomorphic to the third Stasheff associahedron and this remains unaltered in the quantum setting.  In Figure~\ref{hexagons} we replace each cluster by the triangulation of a hexagon determined by the minor indices, as before.

We see that no localisation of coefficients is required and since scalar multiples of all generators of $\qea{n^{-}_{\{\mathrm{1,2}\}}}$ occur in the list of cluster variables, we conclude that $\qea{n^{-}_{\{\mathrm{1,2}\}}}$ is a quantum cluster algebra.

\pagebreak

\bibliographystyle{elsarticle-num}
\bibliography{references}\label{references}

\def\cprime{$'$} \newcommand{\noopsort}[1]{}
\begin{thebibliography}{10}
\expandafter\ifx\csname url\endcsname\relax
  \def\url#1{\texttt{#1}}\fi
\expandafter\ifx\csname urlprefix\endcsname\relax\def\urlprefix{URL }\fi
\expandafter\ifx\csname href\endcsname\relax
  \def\href#1#2{#2} \def\path#1{#1}\fi

\bibitem{FZ-CA1}
S.~Fomin, A.~Zelevinsky, Cluster algebras. {I}. {F}oundations, J. Amer. Math.
  Soc. 15~(2) (2002) 497--529 (electronic).

\bibitem{FZ-CA2}
S.~Fomin, A.~Zelevinsky, Cluster algebras. {II}. {F}inite type classification,
  Invent. Math. 154~(1) (2003) 63--121.

\bibitem{BFZ-CA3}
A.~Berenstein, S.~Fomin, A.~Zelevinsky, Cluster algebras. {III}. {U}pper bounds
  and double {B}ruhat cells, Duke Math. J. 126~(1) (2005) 1--52.

\bibitem{FZ-CA4}
S.~Fomin, A.~Zelevinsky, Cluster algebras. {IV}. {C}oefficients, Compos. Math.
  143~(1) (2007) 112--164.

\bibitem{Fomin-Reading}
S.~Fomin, N.~Reading, Root systems and generalized associahedra, in: Geometric
  combinatorics, Vol.~13 of IAS/Park City Math. Ser., Amer. Math. Soc.,
  Providence, RI, 2007, pp. 63--131.

\bibitem{Keller-CAsurvey}
B.~Keller, Cluster algebras, quiver representations and triangulated
  categories, preprint (2008).
\newblock \href {http://arxiv.org/abs/0807.1960} {\path{arXiv:0807.1960}}.

\bibitem{Caldero-Chapoton}
P.~Caldero, F.~Chapoton, Cluster algebras as {H}all algebras of quiver
  representations, Comment. Math. Helv. 81~(3) (2006) 595--616.

\bibitem{BMRRT}
A.~B. Buan, R.~Marsh, M.~Reineke, I.~Reiten, G.~Todorov, Tilting theory and
  cluster combinatorics, Adv. Math. 204~(2) (2006) 572--618.

\bibitem{GLS-PFV}
C.~Gei\ss, B.~Leclerc, J.~Schr{\"o}er, Partial flag varieties and preprojective
  algebras, Ann. Inst. Fourier (Grenoble) 58~(3) (2008) 825--876.

\bibitem{Hernandez-Leclerc}
D.~Hernandez, B.~Leclerc, Cluster algebras and quantum affine algebras,
  preprint (2009).
\newblock \href {http://arxiv.org/abs/0903.1452} {\path{arXiv:0903.1452}}.

\bibitem{Scott-Grassmannians}
J.~S. Scott, Grassmannians and cluster algebras, Proc. London Math. Soc. (3)
  92~(2) (2006) 345--380.

\bibitem{BZ-QCA}
A.~Berenstein, A.~Zelevinsky, Quantum cluster algebras, Adv. Math. 195~(2)
  (2005) 405--455.

\bibitem{BraidedEnvAlgs}
J.~E. Grabowski, Braided enveloping algebras associated to quantum parabolic
  subalgebras, preprint (2007).
\newblock \href {http://arxiv.org/abs/0706.0455} {\path{arXiv:0706.0455}}.

\bibitem{GAP4}
The GAP~Group, {GAP -- Groups, Algorithms, and Programming, Version 4.4},
  \mbox{\url{http://www.gap-system.org/}}.

\bibitem{QuaGroup}
Willem A. de~Graaf, QuaGroup, a package for GAP4, version 1.3. Software web
  homepage: \mbox{\url{http://www.science.unitn.it/~degraaf/quagroup.html}}.

\bibitem{Gr2nSchubertQCA}
J.~E. Grabowski, S.~Launois, Quantum cluster algebra structures on quantum
  {G}rassmannians and their quantum {S}chubert cells: the finite-type cases,
  preprint (2009).
\newblock \href {http://arxiv.org/abs/0912.4397} {\path{arXiv:0912.4397}}.

\bibitem{Nakajima-QuiverVarClusterAlgs}
H.~Nakajima, Quiver varieties and cluster algebras, preprint (2009).
\newblock \href {http://arxiv.org/abs/0905.0002} {\path{arXiv:0905.0002}}.

\bibitem{Lampe}
P.~Lampe, A quantum cluster algebra of {K}ronecker type and the dual canonical
  basis, preprint (2010).
\newblock \href {http://arxiv.org/abs/1002.2762} {\path{arXiv:1002.2762}}.

\bibitem{LusztigBook}
G.~Lusztig, Introduction to quantum groups, Vol. 110 of Progress in
  Mathematics, Birkh\"auser Boston Inc., Boston, MA, 1993.

\bibitem{Borel}
A.~Borel, Linear algebraic groups, 2nd Edition, Vol. 126 of Graduate Texts in
  Mathematics, Springer-Verlag, New York, 1991.

\bibitem{Brown-Goodearl}
K.~A. Brown, K.~R. Goodearl, Lectures on algebraic quantum groups, Advanced
  Courses in Mathematics. CRM Barcelona, Birkh\"auser Verlag, Basel, 2002.

\bibitem{Kolb}
S.~Kolb, The {AS}-{C}ohen-{M}acaulay property for quantum flag manifolds of
  minuscule weight, J. Algebra 319~(8) (2008) 3518--3534.

\bibitem{Jantzen}
J.~C. Jantzen, Lectures on quantum groups, Vol.~6 of Graduate Studies in
  Mathematics, American Mathematical Society, Providence, RI, 1996.

\bibitem{Leclerc-Zelevinsky}
B.~Leclerc, A.~Zelevinsky, Quasicommuting families of quantum {P}l\"ucker
  coordinates, in: Kirillov's seminar on representation theory, Vol. 181 of
  Amer. Math. Soc. Transl. Ser. 2, Amer. Math. Soc., Providence, RI, 1998, pp.
  85--108.

\bibitem{Scott-QMinors}
J.~Scott, Quasi-commuting families of quantum minors, J. Algebra 290~(1) (2005)
  204--220.

\bibitem{KLR}
A.~C. Kelly, T.~H. Lenagan, L.~Rigal,
  \href{http://dx.doi.org/10.1142/S0219498804000630}{Ring theoretic properties
  of quantum {G}rassmannians}, J. Algebra Appl. 3~(1) (2004) 9--30.
\newblock \href {http://dx.doi.org/10.1142/S0219498804000630}
  {\path{doi:10.1142/S0219498804000630}}.
\newline\urlprefix\url{http://dx.doi.org/10.1142/S0219498804000630}

\bibitem{Magma}
W.~Bosma, J.~Cannon, C.~Playoust, The {M}agma algebra system. {I}. {T}he user
  language, J. Symbolic Comput. 24~(3-4) (1997) 235--265, computational algebra
  and number theory (London, 1993).

\bibitem{Heckenberger-Kolb}
I.~Heckenberger, S.~Kolb, The locally finite part of the dual coalgebra of
  quantized irreducible flag manifolds, Proc. London Math. Soc. (3) 89~(2)
  (2004) 457--484.

\end{thebibliography}

\end{document}